\documentclass[11pt,a4paper,twoside]{amsart}
\usepackage{amssymb,amsmath,amsthm}
\theoremstyle{plain}
\newtheorem{theorem}{Theorem}[section]

\newtheorem{lemma}[theorem]{Lemma}

\newtheorem{proposition}[theorem]{Proposition}
\newtheorem{assumption}[theorem]{Assumption}
\theoremstyle{remark}
\newtheorem{remark}[theorem]{Remark}

\usepackage[dvipdfm,
   bookmarks=true,
   bookmarksnumbered=false,
   bookmarkstype=toc]
{hyperref}

\numberwithin{equation}{section}

\newcommand{\C}{\mathbb{C}}
\newcommand{\R}{\mathbb{R}}

\newcommand{\N}{\mathbb{N}}

\newcommand{\F}{\mathcal{F}}

\renewcommand{\Re}{\operatorname{Re}}
\newcommand{\I}{\infty}
\newcommand{\abs}[1]{\left\lvert #1\right\rvert}
\newcommand{\norm}[1]{\left\lVert #1\right\rVert}
\newcommand{\Lebn}[2]{\left\lVert #1 \right\rVert_{L^{#2}}}
\newcommand{\Sobn}[2]{\left\lVert #1 \right\rVert_{H^{#2}}}

\newcommand{\Jbr}[1]{\left\langle #1 \right\rangle}

\newcommand{\tildeE}[1]{\widetilde{#1}^\varepsilon}

\newcommand{\ceil}[1]{\left\lceil#1\right\rceil}
\newcommand{\IN}{\quad\text{in }}

\def\({\left(}
\def\){\right)}
\def\<{\left\langle}
\def\>{\right\rangle}
\def\le{\leqslant}
\def\ge{\geqslant}

\def\d{{\partial}}
\def\l{\lambda}

\newcommand{\g}{\gamma}

\newcommand{\eps}{\varepsilon}

\begin{document}
\title[WKB approximation for Schr\"odinger-Poisson system]{
Large time WKB approximation for multi-dimensional semiclassical
Schr\"odinger-Poisson system}
\author[S. Masaki]{Satoshi Masaki}
\address{Division of Mathematics\\
Graduate School of Information Sciences\\
Tohoku University\\
Sendai 980-8579, Japan}
\email{masaki@ims.is.tohoku.ac.jp}
\begin{abstract}
We consider the semiclassical Schr\"odinger-Poisson system with 
a special initial data of WKB type such that the solution of 
the limiting hydrodynamical equation
becomes time-global in dimensions at least three.
We give an example of such initial data in the focusing case
via the analysis of the compressible Euler-Poisson equations.
This example is a large data with radial symmetry,
and is beyond the reach of the previous results because the phase part decays
too slowly. 
Extending previous results in this direction, we justify
the WKB approximation of the solution 
with this data for an arbitrarily large interval of $\R_+$.
\end{abstract}
\maketitle
\section{introduction}
This paper is devoted to the study of the semiclassical limit $\eps \to 0$
for the Cauchy problem of the semiclassical Schr\"odinger-Poisson system
for $(t,x) \in \R_+ \times \R^n$
\begin{equation}\label{eq:r3}
\left\{
\begin{aligned}
	&i\eps\d_t u^\eps +\frac{\eps^2}{2}\Delta u^\eps =\l V^\eps_{\mathrm{P}} u^\eps, \\
	&-\Delta V^\eps_{\mathrm{P}} = |u^\eps|^2, \quad V^\eps_{\mathrm{P}}\in L^\I(\R^n), \quad V^\eps_{\mathrm{P}} \to 0 \text{ as }|x|\to\I, \\
	& {u^\eps}_{|t=0}(x) = A_0^\eps (x)e^{i\frac{\Phi_0(x)}{\eps}},
\end{aligned}
\right.
\end{equation}
where $n\ge 3$, $\eps$ is a positive parameter which corresponds to the scaled Planck constant, and $\l$ is a real number.
In addition, the ``initial amplitude'' $A_0^\eps$ is complex-valued and 
the ``initial phase'' $\Phi_0$ is real-valued.
Precise assumption on them is in Assumption \ref{asmp:main}. 
It is known that, if $n\ge 3$ then $V^\eps_{\mathrm{P}}\in L^\I(\R^n)$ is uniquely determined
from $u^\eps \in L^2 \cap L^\I$ as
\[
	c_n (|x|^{n-2}*|u^\eps|^2),
\]
where $c_n$ is a positive constant.
Therefore, the Schr\"odinger-Poisson system \eqref{eq:r3} can be regarded
as a special case of the Hartree equation
\begin{equation}\label{eq:HE}
	i\eps\d_t u^\eps +\frac{\eps^2}{2}\Delta u^\eps =\l (|x|^{\g}*|u^\eps|^2) u^\eps.
\end{equation}
For the well-posedness results on \eqref{eq:r3} and \eqref{eq:HE}
for fixed $\eps>0$, see \cite{CazBook} and references therein.

In this paper, we are interested in the WKB type approximation
for the solution of \eqref{eq:r3}:
\begin{equation}\label{eq:WKB}
	u^\eps(t,x) \sim e^{i\frac{\phi(t,x)}{\eps}}(\beta_0(t,x)
	+ \eps \beta_1(t,x) + \eps^2 \beta_2 + \cdots)
\end{equation}
as $\eps\to0$. 
One way to justify \eqref{eq:WKB} is to employ a modified Madelung transform
\begin{equation}\label{eq:mMT}
	u^\eps(t,x) = a^\eps(t,x) \exp\(i\frac{\phi^\eps(t,x)}{\eps}\)
\end{equation}
and consider the system
\begin{equation}\label{eq:QSP}
\left\{
	\begin{aligned}
	&\partial_t a^\eps + \nabla \phi^\eps \cdot \nabla a^\eps + \frac12 a^\eps \Delta \phi^\eps =i\frac{\eps}{2} \Delta a^\eps, 
	&& a^\eps(0,x)= A^\eps_0(x); \\
	&\partial_t \phi^\eps +\frac12 |\nabla \phi^\eps|^2 + \l V^\eps_{\mathrm{P}} = 0 ,
	&& \phi^\eps(0,x) = \Phi_0(x); \\
	& -\Delta V^\eps_{\mathrm{P}} = |a^\eps|^2, \quad V^\eps_{\mathrm{P}}\in L^\I(\R^n),
	&&V^\eps_{\mathrm{P}} \to 0 \text{ as }|x|\to\I.	
	\end{aligned}
\right.
\end{equation}
Note that $a^\eps$ takes complex value.
Our strategy is to obtain an expansion like
\[
	a^h = a_0 + \eps a_1 + \eps^2 a_2 + \cdots, \quad
	\phi^h = \phi_0 + \eps \phi_1 + \eps^2 \phi_2 + \cdots,
\]
which yields \eqref{eq:WKB} together with \eqref{eq:mMT}.
This method is first applied to analytic data \cite{PGEP}
and to Sobolev data \cite{Grenier98} for certain class of defocusing nonlinearities,
and is generalized to other local nonlinearities in 
\cite{AC-BKW,AC-GP,CR-CMP} and to some nonlocal nonlinearities in
\cite{AC-SP,CM-AA,LL-EJDE,LT-MAA}.
For this method, see also \cite{CaBook,GLM-TJM}.
One verifies that the principal part $(a_0,\phi_0)$ of $(a^\eps,\phi^\eps)$ solves, at least formally,
\begin{equation}\label{eq:limsys2}
\left\{
  \begin{aligned}
    &\partial_t a_0 + \nabla\phi_0 \cdot \nabla a_0 + \dfrac{1}{2}a_0 \Delta \phi_0=0,&&
	a_{0|t=0} = A_0, \\ 
	&\partial_t \phi_0 + \frac12 |\nabla \phi_0|^2 + \l V_{\mathrm{P}}  = 0, && \phi_{0|t=0}=\Phi_0, \\
    &-\Delta V_{\mathrm{P}} = |a_0|^2, \quad V_{\mathrm{P}}\in L^\I(\R^n),
    &&V_{\mathrm{P}} \to 0 \text{ as }|x|\to\I,
  \end{aligned}
\right.
\end{equation}
where $A_0:=\lim_{\eps\to0}A^\eps_0$.
In general, the classical solution of \eqref{eq:limsys2} breaks down in finite time
by a formation of singularity.
The space-time set where the solution ceases to be smooth is called caustic.
At the caustic, the WKB type approximation \eqref{eq:WKB} also breaks down.
The shape of the caustic set depends on the initial data of \eqref{eq:limsys2}.
The aim of this paper is to justify the large time WKB  with 
a special initial data of WKB type which does not cause the caustic.

Whether the caustic phenomena occurs or not 
boils down to the problem of global existence of the classical solution to \eqref{eq:limsys2}.
By choosing $\rho:=|a_0|^2$ and $v:=\nabla \phi_0$, we find
that \eqref{eq:limsys2} becomes the compressible Euler-Poisson equations (see \eqref{eq:EP}, below).
The classical solutions of the compressible Euler-Poisson equations
are studied in \cite{Chae0803,CT-CMS,ELT-IUMJ,MaRIMS}.
We see from \cite{ELT-IUMJ} that there is an example of the initial data which does
not cause the caustic, provided $n=1$ and $\l>0$ (repulsive case, or defocusing case).
For such initial data, the large time WKB analysis of \eqref{eq:r3} 
is shown in \cite{LT-MAA}.
It is pointed out in \cite{MaRIMS} that, for $n\ge3$ and under certain conditions such as radial symmetry, such example exists if
$\l<0$ (attractive case, or focusing case).
Our results are based on this respect.

\subsection{Main result}
We denote by $H^s(\R^n)$ the usual Sobolve space: $H^s(\R^n)
=\{f\in L^2(\R^n); (1-\Delta)^{s/2}f \in L^2(\R^n) \}$.
Let us write $H^s=H^s(\R^n)$, for short.
\begin{assumption}\label{asmp:main}
Suppose $n\ge 3$ and $\l<0$ (focusing case).
Let ``expansion level'' $N$ be a positive integer.
We suppose the following conditions with some $s>n/2+2N+1$:
\begin{enumerate}
\item The initial amplitude $A_0^\varepsilon \in H^{s+1}(\R^n)$ writes
\[
	A_0^\eps = A_0 + \sum_{j=1}^{N} \eps^j A_j + O(\eps^{N+1}) \IN H^{s+1}(\R^n).
\]
Namely, there exist $A_j \in H^{s+1}$ ($j=1,\dots,N$) such that
$\limsup_{\eps\to0}\|A_0^\eps-\sum_{j=0}^N \eps^j A_j\|_{H^{s+1}}/\eps^{N+1}<\I$.
\item The limit $A_0$ of the initial amplitude is radially symmetric, that is, $A_0(x)=A_0(|x|)$: $\R_+\to \C$.
Moreover, there exist $\kappa \ge \ceil{s}+3 -n/2$ and $\delta \in (0,1/4]$ such that
$A_0\in C^{\ceil{s}+3}((0,\I))$ satisfies
\begin{equation}\label{asmp:A7}
	\begin{aligned}
	& A_0(r)\neq 0, \quad \forall r>0;  &
	&\limsup_{r\to0}\frac{|A_0^{(j)}(r)|}{r^{\kappa-j}} <\I 
	\end{aligned}
\end{equation}
for all $j \in [0,\ceil{s}+3]$, and that
the following limits exist and are nonzero:
\begin{equation}\label{asmp:A8}
	 \lim_{r\to0} \frac{A_0(r)}{r^{\kappa}}, \quad
	 \lim_{r\to\I} r^{\frac{n}2+\ceil{s}+3+\delta} A_0^{(\ceil{s+3})}(r),
\end{equation}
where $\ceil{s}$ denotes the minimum integer larger than or equal to $s$. 
\item The initial phase
$\Phi_0$ is a radial function $\Phi_0(x)=\Phi_0(|x|)$: $\R_+ \to \R$ 
given from $A_0$ by the formula
\begin{equation}\label{eq:defgphi0}
	\Phi_0(r) = \int_0^{r}\sqrt{\frac{2|\l|}{(n-2)s^{n-2}}\int_0^{s}|A_0(\sigma)|^2\sigma^{n-1}d\sigma}ds
 +\mathrm{const}.
\end{equation}
\end{enumerate}
\end{assumption}
We now state our main result of this paper.
\begin{theorem}\label{thm:main}
Let Assumptions \ref{asmp:main} be satisfied.
Then, for any $T>0$ and $\eps>0$ with $\eps \le C_1 e^{-C_2  T}$,
there exist a solution $u^\eps \in C([0,T]; H^{s+1})$ of \eqref{eq:r3},
and there also exist $\phi_0 \in C([0,\I);C^{2N+6})$ and
$\beta_j \in C([0,\I); H^{s-2j+3})$ such that the approximation
\begin{equation}\label{eq:WKBHs2}
	u^\eps = e^{i\frac{\phi_0}{\eps}} (\beta_0 + \eps \beta_1
	+\cdots + \eps^{N-1} \beta_{N-1} + O(\eps^{N})) \IN
	L^\I([0,T]; H^{s-2N+1})
\end{equation}
holds, where $C_1$ and $C_2$ are positive constants.
Moreover, $C_2$ is independent of $N$.
\end{theorem}
In other words, this theorem tells us that,
for any $0<\eps\ll1$, the solution $u^\eps$ of \eqref{eq:r3} exists and
\eqref{eq:WKBHs2} is valid for
$t\le C_2^{\prime} \log \frac1\eps +C_1^\prime$.

\begin{remark}
We remark that the initial amplitude $A_0^\eps$ is not necessarily a radial function,
though its limit $A_0$ is supposed to be radially symmetric.
A simple example of $A_0$ which satisfies Assumption \ref{asmp:main} is
	$A_0(r) = r^{\kappa} \psi(r) + r^{-\frac{n}2-\delta}(1-\psi(r))$,
where $\psi(r) \in C^\I((0,\I))$ is a function such that
$0\le \psi \le 1$, $\psi(r)=1$ for $r<1$, and $\psi(r)=0$ for $r>2$.
\end{remark}
\begin{remark}
The initial data which satisfies Assumption \ref{asmp:main} is not necessarily
a small data. Indeed, if some $(A_0.\Phi_0)$ satisfies this assumption and
if $\alpha$ is a complex number,
then $(\alpha A_0, |\alpha| \Phi_0)$ also satisfies this assumption.
\end{remark}
\begin{remark}\label{rmk:samplePhi}
Notice that $\Phi_0$ given by \eqref{eq:defgphi0} belongs to
$L^\I(\R^n)$ if and only if $n\ge5$.
One also sees that
$\nabla \Phi_0$ belongs to $L^p(\R^n)$ only if $p>2^*$.
This is due to the lack of decay at spatial infinity.
\end{remark}

Theorem \ref{thm:main} follows from Theorems \ref{thm:global} and \ref{thm:WKB},
below.
To state them, we make some definitions and notation.
For $n \ge 3$, $s>n/2+1$, $p\in [1, \I]$, and $q \in[1,\I] $,
we define a function space $Y^s_{p,q}(\R^n)$ by
\begin{equation}\label{def:Y}
	Y^s_{p,q}(\R^n) = \overline{C_0^\I (\R^n)}^{\norm{\cdot}_{Y^s_{p,q}(\R^n)}}
\end{equation}
with norm
\begin{equation}\label{def:Y2}
	\norm{\cdot}_{Y^s_{p,q}(\R^n)} :=
	\norm{\cdot}_{L^p(\R^n)} + \norm{\nabla \cdot}_{L^q(\R^n)}
	+ \norm{\nabla^2 \cdot}_{H^{s-2}(\R^n)}.
\end{equation}
We denote $Y^s_{p,q}=Y^s_{p,q}(\R^n)$, for short.
This space $Y^s_{p,q}$, introduced in \cite{MaCascade}, is
a modification of the Zhidkov space $X^s$,
which is defined, for $s>n/2$, by $X^s(\R^n) := \{  f \in L^\I(\R^n) | \nabla f \in H^{s-1}(\R^n) \}$.
The Zhidkov space was introduced in \cite{Zhidkov} (see, also \cite{Gallo}).
Roughly speaking, the exponents $p$ and $q$ in $Y^s_{p,q}$ indicate the decay rates at spatial infinity
of the function and of its first derivative, respectively.
Moreover, $X^s \fallingdotseq Y^s_{\I,2}$ if $n \ge 3$ and
$Y^s_{2,2}=H^s$.
We use the following notation:
\begin{align*}
	Y^s_{I_1,q} &{}:= \cap_{p^\prime\in I_1} Y^s_{p^\prime,q}, &
	Y^s_{p,I_2} &{}:= \cap_{q^\prime\in I_2} Y^s_{p,q^\prime},
\end{align*}
where $I_1$ and $I_2$ are intervals of $[ 1,\I]$.
These notation are sometimes used simultaneously. For example,
$Y^s_{I_1,I_2} := \cap_{p^\prime\in I_1,q^\prime\in I_2} Y^s_{p^\prime,q^\prime}$.
For $q\in [1,n)$, we use the notation $q^* = nq/(n-q)$. The following two relations are 
sometimes useful:
$q<q^*<\I$; $q_1^* > q_2^*$ if and only if $q_1 > q_2$.
If $f \in Y^s_{p,q}$ 
then $|\nabla f| \to 0$ as $|x|\to\I$ by definition,
and so
$\Lebn{\nabla f }{\I}\le \Sobn{\nabla^2 f}{s}<\I$ by Sobolev embedding,
which means $Y^s_{p,q}=Y^s_{p,[q,\I]}$.
Similarly, it follows from the Sobolev embedding and Lemma \ref{lem:Zhidkov1},
below, that
\begin{align*}
	Y^s_{p,q} &{}= Y^s_{p,[\min(q,2^*),\I]},  &
	Y^s_{p,q} &{}\subset Y^s_{q^*,q} \quad\text{ if } q<n.
\end{align*}

We first claim that \eqref{eq:limsys2} has a radial global solution under
Assumption \ref{asmp:main}.
For a solution $(a_0,\phi_0)$ of \eqref{eq:limsys2}, we here introduce
\begin{equation}\label{def:eta}
	\eta(T)=\eta(T;s,p,q)= \sup_{t\in [0,T]}\(\norm{a_0(t)}_{H^{s+3}(\R^n)}
	+\norm{\nabla \phi_0(t)}_{Y^{s+4}_{p,q}(\R^n)}\),
\end{equation}
which is the key value for combining the following two theorems.
\begin{theorem}\label{thm:global}
Let Assumption \ref{asmp:main} be satisfied,
Then, there exists a radial global solution
$(a_0,\phi_0)(t,x)=(a_0, \phi_0)(t,|x|)$ of \eqref{eq:limsys2} satisfying
$a_0(t) \in H^{s+3}(\R^n)$, $\phi_0(t)\in C^{2N+6}$,
and $\nabla \phi_0(t) \in Y^{s+4}_{(2^*,\I],(2,\I]}$
for all $t\ge0$. 
Moreover, for any $p_0>2^*$ and $q_0>2+{4\delta}/{n}$,
it holds that
 \begin{align}
 	&\norm{a_0(t)}_{L^{2}(\R^n)} = \norm{A_0}_{L^{2}(\R^n)},  \quad
 	\norm{\nabla a_0(t)}_{H^{s+2}(\R^n)} = o(1), \label{eq:ordera}\\
	&\norm{\nabla \phi_0(t)}_{Y^{s+4}_{[p_0,\I].[q_0,\I]} (\R^n)} {}= o(1)
	\label{eq:orderphi}
\end{align}
as $t\to\I$. In particular, $\eta(T;s,p,q)=O(1)$ as $T\to\I$
for all $p>2^*$ and $q>2+4\delta/n$.
\end{theorem}
\begin{remark}
Under the following three assumptions;
$($i$)$ $n\ge 3$ and the radial symmetry;
$($ii$)$ $A_0 \in L^2(\R^n)$;
$($iii$)$ $\nabla \Phi(0)=0$ (by the radial symmetry) and $\nabla \Phi_0$ decreases at spatial infinity;
the classical solution $(a_0,\phi_0)$ of \eqref{eq:limsys2}
is time-global if and only if $\l<0$ and $\Phi_0$ is given by \eqref{eq:defgphi0}
(see Theorem \ref{thm:globalCT}).
\end{remark} 

We now consider the WKB analysis of \eqref{eq:r3}.
Let us go back to the equation \eqref{eq:QSP}.
It is proven in \cite{AC-SP} that this system has a local solution 
$(a^\eps,\phi^\eps)$ for $0\le \eps \ll1$
and the solution can be expanded as
\[
	a^h = a_0 + \eps a_1 + \eps^2 a_2 + \cdots, \quad
	\phi^h = \phi_0 + \eps \phi_1 + \eps^2 \phi_2 + \cdots
\]
for a class of initial data of WKB type (see, also \cite{CM-AA}).
However, the initial data satisfying Assumption \ref{asmp:main} is 
out of framework of this results because
the spatial decay of the phase function $\Phi_0$ is slow.
So, we extend the result in this direction (see Remark \ref{rmk:decayPhi}, below).

\begin{assumption}\label{asmp:existence}
Let $n \ge 3$ and $\l \in \R$. Let $N$ be a positive integer.
We suppose the following conditions with some $s>n/2+2N+1$:\\
$\bullet$ The initial amplitude $A_0^\varepsilon$ 
satisfies $(1)$ of Assumption \ref{asmp:main}. \\
$\bullet$  The initial phase $\Phi_0 \in C^{2N+4} (\R^n)$ satisfies
$\nabla \Phi_0 \in Y^{s+2}_{p,q}(\R^n)$ for some $p\in (2^*,\I]$
and $q\in(2,n)$ with $p \ge q$.
\end{assumption}
\begin{remark}
Assumption \ref{asmp:existence} is weaker than Assumption \ref{asmp:main}.
In particular, $A_0$ and $\Phi_0$ are not necessarily radially symmetric.
\end{remark}
\begin{theorem}\label{thm:WKB}
Let Assumption \ref{asmp:existence} be satisfied.
Suppose that \eqref{eq:limsys2} has a global solution $(a_0,\phi_0)$
which satisfies $\eta(T;s,p,q)<\I$ for all $T<\I$.
Then, for any $T>0$ and $\eps>0$ with $\eps \le C_1 \eta(T) e^{-C_2\eta(T)T}$, 
there exits a solution $u^\eps \in C([0,T]; H^{s+1})$, and
there also exist $\beta_j\in C([0,\I); H^{s-2j+1})$ such that
the approximation
\begin{equation}\label{eq:WKBHs}
	u^\eps = e^{i\frac{\phi_0}{\eps}} (\beta_0 + \eps \beta_1
	+\cdots + \eps^{N-1} \beta_{N-1} + O(\eps^{N})) \IN
	L^\I([0,T]; H^{s-2N+1})
\end{equation}
holds, where $C_1$ and $C_2$ are positive constants.
Moreover, $C_2$ depends only on $n$ and $s$.
\end{theorem}
\begin{remark}\label{rmk:decayPhi}
In \cite{AC-SP} and \cite{CM-AA},
the assumption on the initial phase $\Phi_0$ is that
$\Phi_0 \in L^\I$, 
$\nabla \Phi_0 \in L^{2^*}\cap L^{n-}\cap L^\I$, and $\nabla^2 \Phi_0 \in L^2\cap L^\I$ (up to a subquadratic polynomial).
Here, we assume
$\nabla \Phi_0 \in L^{2^*+}\cap L^\I$ and $\nabla^2 \Phi_0 \in L^{2+}\cap L^\I$ at best.
Remark that in our framework, $\Phi_0$ does not necessarily 
belong to any Lebesgue space nor is not a polynomial, as in \eqref{eq:defgphi0}.
\end{remark}
\begin{remark}\label{rmk:expansionlevel}
In Theorem \ref{thm:WKB}, the constant $C_2$ is independent of the 
$N$ and the initial data.
This point is an improvement because
this constant was proportional to $N$ in \cite{LT-MAA}.
On the other hand, $C_1$ depends on them.
However, the influence of $C_1$ is much smaller than of $C_2$.
\end{remark}
\smallbreak

The rest of this paper is organized as follows:
We prove Theorems \ref{thm:global} and \ref{thm:WKB}
in Sections \ref{sec:proof2} and \ref{sec:proof1}, respectively.
Sections \ref{sec:pre2} and \ref{sec:pre1} are devoted to preliminary results for the proofs.
In Appendix \ref{sec:EP}, we prove some results on the radial compressible
Euler-Poisson equations which we use for the proof of Theorem \ref{thm:global}.
\section{Preliminaries for the proof of Theorem \ref{thm:global}}\label{sec:pre2}
In this section, we collect some preliminary results which will be used for the proof of 
Theorem \ref{thm:global}.
In Subsections \ref{subsec:EP} and \ref{subsec:limsys2},
we give an explicit example of the global solution to \eqref{eq:limsys2}
based on results in \cite{MaRIMS}. 
We also show an elementary equality in Subsection \ref{subsec:equality}.

\subsection{Global existence of the solution to the compressible Euler-Poisson system}\label{subsec:EP}
Let $(a_0,\phi_0)$ be a solution of \eqref{eq:limsys2}.
Then, one easily sees that $\rho:=|a_0|^2$ and $v:=\nabla \phi_0$ solve
the compressible Euler-Poisson equations
\begin{equation}\label{eq:EP}
\left\{
  \begin{aligned}
    &\partial_t \rho + \mathrm{div} (\rho v)=0,&&
	\rho_{|t=0} = |A_0|^2, \\ 
	&\partial_t v + (v\cdot \nabla) v + \l \nabla V_{\mathrm{P}}  = 0, && v_{|t=0}=\nabla \Phi_0, \\
    &-\Delta V_{\mathrm{P}} = \rho, \quad V_{\mathrm{P}}\in L^\I(\R^n),
    &&V_{\mathrm{P}} \to 0 \text{ as }|x|\to\I.
  \end{aligned}
\right.
\end{equation}
Now, we assume the radial symmetry and consider
the radial version of \eqref{eq:EP} 
\begin{equation}\label{eq:rEP}
\left\{
\begin{aligned}
	&\rho_t + r^{-(n-1)}\d_r (r^{n-1} \rho v) = 0, 
	&&\rho(0,r) = \rho_0(r):=|A_0(r)|^2;\\
	&v_t + v \d_r v + \l \d_r V_{\mathrm{P}} =0,
	&&v(0,r) = v_0(r):=\Phi_0^\prime(r);\\
	&-r^{-(n-1)}\d_r (r^{n-1} V_{\mathrm{P}}) = \rho,\, V_{\mathrm{P}} \in L^\I, && V_{\mathrm{P}}\to 0\text{ as }r\to \I,
\end{aligned}
\right.
\end{equation}
where unknowns are now real-valued functions
$\rho: \R_+^2 \to \R_+$ and $v:\R_+^2\to \R$.
Let us introduce several function spaces.
For a nonnegative integer $k$, we define
\begin{equation}\label{def:Dk}
	D^k :=
	\begin{cases}
	C([0,\I))& \text{ if } k=0, \\
	C([0,\I))\cap C^{k}((0,\I)) & \text{ if } k\ge1.
	\end{cases}
\end{equation}
Similarly, we define
\begin{align*}
	D_{\rho}^k:={} & D^k \cap L^1((0,\I),r^{n-1}dr) , &
	D_{a}^k:={} & D^k \cap L^2((0,\I),r^{n-1}dr)
\end{align*}
for $k \ge 0$ and
\[
	D^k_\phi :=
	\begin{cases}
	C^1([0,\I))& \text{ if } k=1, \\
	C^1([0,\I))\cap C^{k}((0,\I)) & \text{ if } k>1
	\end{cases}
\]
for $k\ge1$.
Let us start with the following theorem announced in \cite{MaRIMS}:
\begin{theorem}[Corollary 1.17 in \cite{MaRIMS}]\label{thm:EP}
Let $n \ge 3$, or $n\ge1$ and $\l<0$.
Suppose $\rho_0\in D^0_\rho$ is not identically zero and
$v_0\in D^1$ satisfies $v_0(0)=0$
and $v_0 \to 0$ as $r\to \I$.
Then, the solution of \eqref{eq:rEP} is global
if and only if  $n\ge 3$, $\l<0$, and the initial data is of particular form
\[
	v_0(r) = \sqrt{\frac{2|\l|}{(n-2) r^{n-2}}\int_0^r \rho_0(s) s^{n-1}ds}.
\]
Suppose $\l<0$ and $n \ge 3$.
If $\rho_0\in D^k_\rho$ for some $k\ge0$ and if $v_0$ is as above,
then $v_0 \in D^{k+1}$ and the corresponding solution is
\begin{align*}
	\rho &{}\in C^2([0,\I),D^k_\rho) \cap C^\I((0,\I),D^k_\rho), \\
	v &{}\in C^1([0,\I),D^{k+1}) \cap C^\I((0,\I),D^{k+1})
\end{align*}
and given explicitly by
\begin{align*}
	\rho(t,X(t,R)) &{}= \rho_0(R) \(1+\dfrac{nv_0(R)}{2R}t\)^{-1}\(1+\dfrac{2|\l| R \rho_0(R)}{(n-2)v_0(R)}t\)^{-1}, \\
	v(t,X(t,R)) &{}= v_0(R)\( 1+ \frac{nv_0(R)}{2R}t \)^{\frac{2}{n}-1},
\end{align*} 
where $X(t,R)=R(1+\frac{nv_0(R)}{2R}t)^{2/n}$.
Moreover, this solution is unique in $C^2([0,\I),D^0) \times C^1([0,\I),D^1)$.
Furthermore, a pair of functions of $(t,x)\in \R_+\times \R^n$ defined as
${\bf r}(t,x):=\rho(t,|x|)$ and ${\bf v}(t,x):=\frac{x}{|x|}v(t,|x|)$
solve \eqref{eq:EP} in the distribution sense.
\end{theorem}
In \cite{MaRIMS}, the proof of this theorem was left incomplete.
We illustrate the proof in Appendix \ref{sec:EP}.
\subsection{Global existence of the solution to \eqref{eq:limsys2}}\label{subsec:limsys2}
We consider the radial version of \eqref{eq:limsys2}:
\begin{equation}\label{eq:rlimsys2}
\left\{
  \begin{aligned}
    &\partial_t a_0 + \partial_r \phi_0 \partial_r a_0 + \dfrac{a_0}{2r^{n-1}} \partial_r( r^{n-1}\partial_r \phi_0) =0,&&
	a_{0|t=0} = A_0, \\ 
	&\partial_t \phi_0 + \frac12 (\d_r\phi_0)^2+ \l V_{\mathrm{P}}  = 0, && \phi_{0|t=0}=\Phi_0, \\
    & -\partial_r(r^{n-1}\partial_r V_{\mathrm{P}}) = r^{n-1}|a_0|^2,
    \quad V_{\mathrm{P}}\in L^\I(\R_+),
    &&V_{\mathrm{P}} \to 0 \text{ as }r\to\I.
  \end{aligned}
\right.
\end{equation}

\begin{theorem}\label{thm:globalCT}
Suppose $n \ge 3$, or $n\ge1$ and $\l<0$.
Suppose $A_0\in D^0_a$ is not identically zero and
$\Phi_0\in D^2_\phi$ satisfies $\Phi_0^\prime(0)=0$ and $\Phi_0^\prime(r) \to 0$ as $r\to \I$.
Then, the solution of \eqref{eq:rlimsys2} is global
if and only if $n\ge 3$, $\l<0$, 
and the initial data is of particular form
\begin{equation}\label{def:globalphi0}
	\Phi_0(r) = \int_0^{r}\sqrt{\frac{2|\l|}{(n-2)s^{n-2}}\int_0^{s}|A_0(\sigma)|^2\sigma^{n-1}d\sigma}ds
 +\mathrm{const}.
\end{equation}
Moreover, if $A_0 \in D^k_a$ for some $k\ge0$, then
the above $\Phi_0$ belongs to $D^{k+2}_\phi$ 
and the corresponding global solution 
\begin{align*}
	a_0 &{}\in C^2([0,\I),D^k_a) \cap C^\I((0,\I),D^k_a) \\
	\phi_0 &{} \in C^1([0,\I),D^{k+2}_\phi) \cap C^\I((0,\I),D^{k+2}_\phi)
\end{align*} 
are given explicitly as
\begin{equation}\label{def:gap}
\begin{aligned}
	a_0(t,X(t,R)) &{}= A_0(R)\(1+\dfrac{nv_0(R)}{2R}t\)^{-\frac12}\(1+\dfrac{2|\l| R |A_0(R)|^2}{(n-2)v_0(R)}t\)^{-\frac12}, \\
	\phi_0(t,X(t,R)) &{}= \Phi_0(R)+ \frac{t}2 \(\Phi_0^\prime(R)^2 + \frac{n-2}{2} \int_0^R \frac{\Phi_0^\prime(r)^2}{r} dr \) + g(t),
\end{aligned} 
\end{equation}
where $X(t,R)=R(1+\frac{n \Phi_0^\prime(R)}{2R}t)^{2/n}$, and $g$ is a function of time given by
\begin{multline*}
	g(t)= \\
	\left\{
	\begin{aligned}
	&\frac{2|\l|}{(n-2)(4-n)}\int_0^\I \frac{|A_0(r)|^2r^2}{\Phi_0^\prime(r)}\left[\(1+\frac{n\Phi_0^\prime(r)}{2r}t\)^{\frac{4}n-1}-1 \right]dr, & \text{ if }n\neq 4,\\
	&|\l| \int_0^\I \frac{|A_0(r)|^2r^2}{\Phi_0^\prime(r)} \log\(1+\frac{2 \Phi_0^\prime(r)}{r}t\)dr, & \text{ if }n=4.
	\end{aligned}
	\right.
\end{multline*}
Furthermore, the solution is unique.
\end{theorem}
This theorem is an immediate consequence of Theorem \ref{thm:EP}
and the following lemma,
which is a modification of \cite[Lemma 3.1]{LT-MAA}.
\begin{lemma}\label{lem:T1=T2}
Let $A_0 \in D^{k}_a$ and $\Phi_0 \in D^{k+2}_\phi$ for some $k\ge0$.
Then, the following two statements are equivalent;
\begin{enumerate}
\item the system \eqref{eq:rlimsys2} has a unique solution $(a_0,\phi_0)$ in $C([0,T),D^k_a\times D^{k+2}_\phi)\cap C^1((0,T),D^k_a\times D^{k+2}_\phi)$ 
with $(a_0,\phi_0)_{|t=0}=(A_0,\Phi_0)$;
\item the radial Euler-Poisson equations \eqref{eq:rEP} has a unique solution
$(\rho,v)$ in $C([0,T),D^k_\rho \times D^{k+1})\cap C^1((0,T),D^k_\rho\times D^{k+1})$ 
with  $(\rho,v)_{|t=0}=(|A_0|^2,\Phi_0^\prime)$.
\end{enumerate}
Moreover, the maximal existence times of $(a_0,\phi_0)$ and of $(\rho,v)$ are the same.
\end{lemma}

\subsection{An equality}\label{subsec:equality}
In the forthcoming section, we will investigate the regularity
of the radial global solution given in Theorem \ref{thm:globalCT}.
Especially, we investigate higher derivatives of $a_0$ and $\phi_0$.
The following equality is useful,
which reflects the special structure of the initial data.
\begin{lemma}
Let $n\ge3$ and $\l<0$. Suppose $A_0 \in C^M((0,\I))$ for large integer $M$.
Define
\[
	v_0(r) = \sqrt{ \frac{2|\l|}{(n-2)r^{n-2}}\int_0^r |A_0(s)|^2s^{n-1} ds}.
\]
Then, there exist real constants $\alpha_{j,k}$ and $\beta_{l,m,k}$ such that
the following equality holds for $k \in[1, M+1]$:
\begin{equation}\label{eq:higherv0}
	\sum_{j=0}^k \alpha_{j,k} r^j v_0^{(j)}
	= \sum_{l=1}^k \sum_{m\in (\N\cup\{0\})^l, |m| \le k-l}
	\beta_{l,m,k} \frac{\(\prod_{i=1}^l \rho_0^{(m_i)}\)r^{2l+|m|}}{v_0^{2l-1}},
\end{equation}
where $\rho_0:=|A_0|^2$, $g^{(m)}$ denotes the $m$-th derivative of $g$ with 
the convention $g^{(0)}=g$.
Moreover, $\alpha_{k,k}=1\neq 0$.
\end{lemma}
\begin{proof}
By definition, $v_0$ satisfies
\begin{equation}\label{eq:1stv0}
	\frac{n-2}2 v_0 + r v_0^\prime =  \frac{|\l|}{n-2} \frac{\rho_0 r^2}{v_0}.
\end{equation}
This implies that \eqref{eq:higherv0} holds if $k=1$ and $a_{1,1}=1$.
Then, operating $r \partial_r$ to the both sides of \eqref{eq:1stv0}
and substituting \eqref{eq:1stv0} in the right hand side,
we obtain \eqref{eq:higherv0} with $k=2$.
Repeating this argument, we obtain the result.
\end{proof} 

\section{Proof of Theorem \ref{thm:global}}\label{sec:proof2}
We are now in a position to prove Theorem \ref{thm:global}.
The global solution has been obtained in previous section.
In Subsection \ref{subsec:initreg}, we check the regularity of the solution 
at the initial time $t=0$.
Then, we investigate the regularity of the solution for $t\ge0$ and 
establish an estimate on the time-order as $t\to\I$ of the norm of the solution
in Subsection \ref{subsec:positivereg}.
\subsection{Regularity at the initial time}\label{subsec:initreg}
Let us first consider the regularity of the initial data of \eqref{eq:limsys2}
which satisfies Assumption \ref{asmp:main}.
\begin{proposition}\label{prop:rega0v0}
Let Assumption \ref{asmp:main} be satisfied.
Then, 
${\bf A}_0(x)=A_0(|x|)$ and
	${\bf \Phi}_0 (x)=\Phi_0 (|x|)$
satisfies ${\bf A}_0 \in H^{s+3}(\R^n)$ and $\nabla {\bf \Phi}_0 \in Y^{s+4}_{(2^*,\I],(2,\I]}(\R^n)$,
respectively.
\end{proposition}
\begin{proof}
{\bf Step 1}.
We first collect the decay property of $A_0$, $\rho_0:=|A_0|^2$,
and $v_0:=\Phi_0^\prime$.
By \eqref{asmp:A7} and \eqref{asmp:A8} and the definition of $v_0$
, there exist positive constants $c$ and $C$ such that
\begin{equation}\label{eq:decayav0}
		\begin{aligned}
		c r^{\kappa} \le{}& |A_0(r)| \le C r^{\kappa}, \\
		c r^{2\kappa} \le{}& \rho_0(r) \le C r^{2\kappa}, \\
		c r^{\kappa+1} \le{}& v_0(r) \le C r^{\kappa+1}
		\end{aligned}
\end{equation}
as $r\to0$.
 Then, we use \eqref{eq:higherv0} with $k=1$ to obtain
\[
	v_0^\prime(r) = \alpha \frac{v_0(r)}{r} + \beta \frac{|A_0(r)|^2r}{v_0(r)}=
	O(r^{\kappa})
\]
as $r\to0$.
Moreover, assumption \eqref{asmp:A7} implies
$r^j \rho_0^{(j)}(r) = O(r^{2\kappa})$ as $r\to0$
for $j \in [1,\ceil{s}+3]$.
By this estimate and \eqref{eq:higherv0} with $k=j$, we see by induction that
\begin{equation}\label{eq:decayv0}
	v_0^{(j)}(r) = O(r^{\kappa+1-j})
\end{equation}
as $r\to0$ holds for $j \in [1,\ceil{s}+4]$. 

We next consider the decay rate as $r\to\I$.
Denote  by $L_\I$ the second limit in \eqref{asmp:A8}.
$L_\I\neq0$ by assumption. It follows by l'H\^opital's rule that
\[
	\lim_{r\to\I} r^{\frac{n}2 + j + \delta} A_0^{(j)}(r)
	=(-1)^{\ceil{s}-j}\frac{\Gamma(\frac{n}{2}+\ceil{s}+3+\delta)}{\Gamma(\frac{n}{2}+j+\delta)}L_\I \neq 0
\]
for $j\in[0,\ceil{s}+3]$, where $\Gamma$ is the Gamma function.
Then, we see that
$\lim_{r\to\I} r^{{n} + j + 2\delta} \rho_0^{(j)}(r)$
exists for $j\in[0,\ceil{s}+3]$. Once they exist, then
\begin{align*}
	\lim_{r\to\I} r^{n+j+2\delta} \rho_0^{(j)}(r)
	{}& = (-1)^j \frac{\Gamma({n}+2\delta)}{\Gamma({n}+j+2\delta)}
	\lim_{r\to\I} r^{n+2\delta} \rho_0(r)\\
	{}&=  (-1)^j \frac{\Gamma({n}+2\delta)}{\Gamma({n}+j+2\delta)}
	\(\frac{\Gamma(\frac{n}{2}+\ceil{s}+3+\delta)}{\Gamma(\frac{n}{2}+\delta)}\)^2 |L_\I|^2 \neq0	
\end{align*}
for all $j \in [0,\ceil{s}+3]$.
In particular, $\rho_0^{(j)}(r)$ is exactly order $r^{-n-j-2\delta}$ as $r\to\I$.
Then, one sees from \eqref{eq:higherv0} that
\[
	\lim_{r\to\I} r^{\frac{n}2-1}\sum_{j=0}^k \alpha_{j,k} r^j v_0^{(j)}(r)
	=0
\]
for all $k\in [1,\ceil{s}+4]$, where $\alpha_{j,k}$ is the same constant as in
\eqref{eq:higherv0}.
Since $\alpha_{k,k}=1\neq 0$, an induction argument proves that
$\lim_{r\to\I} r^{\frac{n}2+j-1}v_0^{(j)}(r)$ ($j=0,1,\dots,\ceil{s}+1$) exists.
Once they exist, they satisfy
\begin{align*}
	\lim_{r\to\I} r^{\frac{n}2+j-1}v_0^{(j)}(r)
	={}&(-1)^{j}\frac{\Gamma(\frac{n}2-1)}{\Gamma(\frac{n}2+j-1)}
	\lim_{r\to\I} r^{\frac{n}2-1}v_0(r) \\
	={}&(-1)^{j}\frac{\Gamma(\frac{n}2-1)}{\Gamma(\frac{n}2+j-1)}
	\sqrt{\frac{2|\l|}{n-2}\int_0^\I|A_0(r)|^2r^{n-1}dr}\neq 0.
\end{align*}
As a result, there exist $r_1$ and constants $c$ and $C$ such that if $r\ge r_1$ then
\begin{equation}\label{eq:decayavinf}
\begin{aligned}
	cr^{-\frac{n}2-j-\delta} \le{}& |A_0^{(j)}(r)| \le Cr^{-\frac{n}2-j-\delta},\\ 
	cr^{-{n}-j-2\delta} \le{}& |\rho_0^{(j)}(r)| \le Cr^{-{n}-j-2\delta}, \\
	cr^{-\frac{n}2+1-j} \le{}& |v_0^{(j)}(r)| \le Cr^{-\frac{n}2+1-j}, \\
\end{aligned}
\end{equation}
for $j \in [0,\ceil{s}+3]$ (the third inequality also holds
for $j=\ceil{s}+4$).
\smallbreak

{\bf Step 2}.
We show ${\bf A}_0 \in H^{s+3}(\R^n)$.
Note that ${\bf A}_0 \in L^{2}(\R^n)$ follows from \eqref{eq:decayavinf}
and \eqref{asmp:A7}.
Moreover, for $k \in [1,\ceil{s}+3]$, we have
\[
	\norm{\nabla^k {\bf A}_0}_{L^2(\R^n)}
	\le C\sum_{j=1}^k \int_0^\I |r^{j-\ceil{s}-3} A_0^{(j)}(r)|^2
	r^{n-1} dr 
	< \I
\]
thanks to \eqref{eq:decayavinf} and \eqref{asmp:A7}.

{\bf Step 3}.
Let us prove that $\nabla {\bf \Phi}_0 \in Y^{s+4}_{(2^*,\I],(2,\I]}(\R^n)$.
It follows from \eqref{eq:decayav0} and \eqref{eq:decayavinf} that
\begin{equation*}
	\norm{\nabla {\bf \Phi}_0}_{L^p(\R^n)}^p =
	\int_0^\I | v_0(r) |^q r^{n-1} dr < \I
\end{equation*}
if $(-n/2+1)p+n-1<-1$, that is, if $p>2^*$.
We next consider $\nabla (\nabla {\bf \Phi}_0)$. Note that
\begin{equation*}
	\norm{ \nabla (\nabla {\bf \Phi}_0)}_{L^q(\R^n)}^q \le C \sum_{j=0}^1
	\int_0^\I \left| r^{j-1} v_0^{(j)}(r) \right|^q r^{n-1} dr.
\end{equation*}
We deduce from \eqref{eq:decayav0}, \eqref{eq:decayv0}, and \eqref{eq:decayavinf} that
the right hand side is bounded if $-(n/2)q+n-1<-1$,
that is, if $q>2$.
The proof of $\nabla^k (\nabla {\bf \Phi}_0)\in L^2(\R^n)$ ($k\in[2,\ceil{s}+4]$)
is similar.
An elementary computation shows that
\[
	\norm{\nabla^k (\nabla {\bf \Phi}_0)}_{L^2(\R^n)}
	\le C \sum_{j=0}^k \int_0^\I \abs{\frac{ v_0^{(j)}(r)}{r^{k-j}}}^2 r^{n-1}dr.
\]
This is finite because of \eqref{eq:decayv0} and \eqref{eq:decayavinf}
and asumpton $\kappa \ge \ceil{s}+3 - {n}/2$.
\end{proof}
\subsection{Persistence of the regularity}\label{subsec:positivereg}
We next show that the global solution given in Theorem \ref{thm:globalCT}
keeps the same regularity as the initial data for all positive time,
thanks to its explicit representation.
\begin{proposition}\label{prop:regav}
Let Assumption \ref{asmp:main} be satisfied.
Let  $(a_0(t,r),\phi_0(t,r))$ be the global solution of \eqref{eq:rlimsys2}
given by the formula \eqref{def:gap} in Theorem \ref{thm:globalCT}.
Then, the corresponding global solution 
\begin{equation}\label{eq:defAvt}
	{\bf a}(t,x) = a_0(t,|x|), \quad
	{\bf \Phi}(t,x) = \phi_0(t,|x|)
\end{equation}
of \eqref{eq:limsys2} satisfies 
	$({\bf a}(t),\nabla{\bf \Phi}(t)) \in H^{s+3}(\R^n)\times Y^{s+4}_{(2^*,\I],(2,\I]}(\R^n)$
for all $t\ge0$.
Moreover, \eqref{eq:ordera} and \eqref{eq:orderphi} hold.
\end{proposition}
\begin{proof}
First of all, we put $v_0(r)=\Phi_0^\prime(r)$ and
\begin{align*}
	F(R) :={}& \frac{nv_0(R)}{2R}\ge0, &
	G(R) :={}& \frac{2|\l| |A_0(R)|^2 R}{(n-2)v_0(R)}=v_0^\prime(R) + \frac{(n-2)v_0(R)}{2R}\ge0.
\end{align*}
Then, they simplify the notation as follows:
\begin{align*}
	a_0(t,X(t,R))={}& A_0(R)( 1+ F(R)t )^{-1/2}( 1+ G(R)t)^{-1/2}, \\
	\d_r \phi_0(t,X(t,R))={}& v_0(R)( 1+ F(R)t  )^{2/n-1}, \\
	X(t,R) ={}& R(1+F(R)t)^{2/n}, \\
	\partial_R X(t,R) ={}& (1+F(R)t)^{2/n-1}(1+G(R)t).
\end{align*}
\smallbreak

{\bf Step 1}.
Let us collect the properties of $F$ and $G$.
First is the decay rate as $R\to0$.
It follows from \eqref{eq:decayav0} and \eqref{eq:decayv0} that
\begin{equation}\label{eq:decayFG0}
	F^{(j)}(R) =O(R^{\kappa-j}), \quad
	G^{(j)}(R)=O(R^{\kappa-j})
\end{equation}
as $R\to0$ for all $j\in[0,\ceil{s}+3]$ (the first estimate also holds for
$j=\ceil{s}+4$).
Notice that \eqref{asmp:A7} gives $G(R)>0$ for all $R>0$.
Therefore, there exist $R_0$ and positive constants $c$ and $C$ such that
\begin{equation}\label{eq:FGratio}
	c \le \frac{F(R)}{R^{\kappa}} \le C, \quad
	c \le \frac{G(R)}{R^{\kappa}} \le C, \quad
	c \le \frac{G(R)}{F(R)} \le C
\end{equation}
holds for $R\le R_0$.
A similar argument as in Step 1 of the proof of Proposition 
\ref{prop:rega0v0} shows that
there exist $R_1$ and positive constants $c$ and $C$ such that
\begin{equation}\label{eq:decayFGinf}
\begin{aligned}
	cR^{-\frac{n}2-j} \le{}& |F^{(j)}(R)| \le CR^{-\frac{n}2-j}, \\
	cR^{-\frac{n}2-j-2\delta} \le{}& |G^{(j)}(R)| \le CR^{-\frac{n}2-j-2\delta} \\
\end{aligned}
\end{equation}
for $j \in [0,\ceil{s}+3]$ if $R\ge R_1$ (the first inequality holds for $j=\ceil{s}+4$).
\smallbreak

{\bf Step 2-a}. We show the uniform boundedness of
$\nabla {\bf \Phi} \in Y^{s+1}_{p,q}(\R^n)$ in time for $p>2^*$
$q>2+4\delta/n$.
Since $\sup_{r\ge0}(|F(r)|+|G(r)|)<\I$, we see that
\begin{align*}
	\norm{\nabla {\bf \Phi}(t)}_{L^p(\R^n)}^p &{}\le
	\int_0^\I |v_0(R)|^p R^{n-1}(1+F(R)t)^{1-\frac{2p}{2^*}}
	(1+G(R)t)dR \\
	&{}\le
	C\norm{ \nabla {\bf \Phi}_0}_{L^p(\R^n)}^p <\I
\end{align*}
for any $t\ge0$ and $p>2^*$, where we have used the fact that
$\frac{1+G(R)t}{1+F(R)t}$ is bounded uniformly in $t$ and $R$.
The $L^\I$ estimate is obtained in the similar way.
By the Lebesgue convergence theorem, one sees that $\norm{\nabla {\bf \Phi}(t)}_{L^p(\R^n)}\to 0$ as $t\to\I$.
\smallbreak

{\bf Step 2-b}.
We next estimate
\begin{align*}
	&\norm{\nabla (\nabla {\bf \Phi})(t)}_{L^q(\R^n)}^q \\
	&{}\le C \int_0^\I \left|\frac{v(t,X)}{X}\right|^qX^{n-1}\d_RX dR 
	+ C \int_0^\I |v^\prime(t,X)|^qX^{n-1}\d_RX dR \\
	&{}\le C \int_0^\I |v_0(R)|^q R^{n-q-1}(1+F(R)t)^{1-q}(1+G(R)t)dR\\
	&\quad{}+C \int_0^\I |v^\prime_0(R)|^qR^{n-1}(1+F(R)t)(1+G(R)t)^{1-q}dR \\
	&\quad{}+ C \int_0^\I |v_0(R)|^q R^{n-1}(1+F(R)t)^{1-q}(1+G(R)t)^{1-q} |F^\prime(R)|^q t^q dR
\end{align*}
for $q>4\delta/n$.
We divide each time integrals as $\int_0^\I=\int_0^{R_0}+\int_{R_0}^{R_1}+\int_{R_1}^\I$. 
We easily handle the seconds and prove
that they are order $o(1)$ as $t\to\I$ since
integrals are over a bounded interval. 
Hence, we left the detail
and only establish uniform estimates of the integrals on the intervals $(0,R_0)$ and $(R_1,\I)$.

By \eqref{eq:FGratio},
$\frac{1+F(R)t}{1+G(R)t}$ is bounded from below and above
uniformly in time for $R \le R_0$.
Therefore,
\begin{align*}
	 &\int_0^{R_0} |v_0(R)|^q R^{n-q-1} (1+F(R)t)^{1-q}
(1+G(R)t)dR \\
	&+\int_0^{R_0} |v^\prime_0(R)|^qR^{n-1}(1+F(R)t)(1+G(R)t)^{1-q}dR \\
	 &{}\le C \int_0^{R_0} R^{\kappa q+n-1}(1+F(R)t)^{2-q}
	 dR.
\end{align*}
Apparently, these two terms are order $o(1)$. Similarly,
$\frac{1+F(R)t}{1+G(R)t}\le C R^{2\delta}$ yields
\begin{align*}
	 &\int_{R_1}^\I |v_0(R)|^q R^{n-q-1} (1+F(R)t)^{1-q}
(1+G(R)t)dR \\
	&+\int_{R_1}^\I |v^\prime_0(R)|^qR^{n-1}(1+F(R)t)(1+G(R)t)^{1-q}dR \\
	 &{}\le C \int_{R_1}^\I R^{-\frac{n}2q+n-1+2\delta}dR <\I
\end{align*}
if $-\frac{n}2q+n+2\delta<0$, that is, if $q>2+4\delta/n$.
On the other hand,
\begin{align*}
	&\int_0^{R_0} |v_0(R)|^q R^{n-1}(1+F(R)t)^{1-q}(1+G(R)t)^{1-q}|F^\prime(R)|^q t^q dR \\
	&{}\le C \int_0^{R_0} R^{\kappa q+n-1} \(R^{\kappa q} (1+F(R)t)^{-q}
	 t^q \)dR.
\end{align*}
It follows from the Young inequality that
\[
	R^{\kappa q} (1+F(R)t)^{2-2q}
	 t^q 
	 \le R^{\kappa q}t^{q} (CR^{\kappa}t)^{-q}
	 \le C.
\]
Therefore, this term is also bounded.
Similarly,
\begin{align*}
	&\int_{R_1}^\I |v_0(R)|^q R^{n-1}(1+F(R)t)^{1-q}(1+G(R)t)^{1-q}|F^\prime(R)|^q t^q dR \\
	&{}\le C \int_{R_0}^\I R^{-\frac{n}2q+n-1+2\delta} \(R^{-\frac{n}2q} (1+F(R)t)^{2-2q}
	 t^q \) (1+G(R)t)^{2-q} dR.
\end{align*}
This is uniformly bounded as $t\to\I$ if $q>2+4\delta/n$
because the following estimate is true for $R\ge R_1$:
\[
	R^{-\frac{n}2q} (1+F(R)t)^{-q}
	 t^q
	 \le C \(\frac{R^{-\frac{n}2}t}{1+R^{-\frac{n}2}t}\)^q .
\]
\smallbreak

{\bf Step 2-c}.
Finally, let us prove $\nabla^k  (\nabla {\bf \Phi}) \in L^2(\R^n)$ for $k\in[2,\ceil{s}+4]$.
Note that
\begin{equation}\label{eq:pf_regv05}
	\norm{\nabla^k  (\nabla {\bf \Phi})(t)}_{L^2(\R^n)}^2 \le C \sum_{j=0}^k 
	\int_0^\I \left| \frac{\partial_r^j v(t,X(t,R))}{X(t,R)^{k-j}}\right|^2
	X^{n-1}\d_RXdR.
\end{equation}
For simplicity, we define $v(t,r):=\d_r \phi(t,r)$.
By the explicit representations of $v_0$ and $\d_R X$, we obtain
\begin{multline}\label{eq:derivativev}
	(\partial_r^k v)(t,X(t,R)) = \sum_{l_i \ge0, l_2<k,l_1+l_2+l_3 \le k} \sum_{m_i\in (\N\cup \{0\})^{l_i},|m_1|+|m_2|= k-l_1-l_2-l_3}\\ C_{k,l_i,m_1,m_2}
	v_0^{(l_3)}(R) (1+F(R)t)^{(k-1)\(1-\frac{2}{n}\)-l_1} 
	(1+G(R)t)^{-k-l_2}t^{l_1+l_2}  \\
	\prod_{i_1=1}^{l_1}F^{(1+m_{1i_1})}(R)
	\prod_{i_2=1}^{l_2}G^{(1+m_{2i_2})}(R).
\end{multline}
Therefore, our task is to prove that
\begin{multline}\label{eq:pf_regv1}
	\int_0^\I \frac{|v_0^{(l_3)}(R)|^2}{ R^{2k-2j-n+1} }
	(1+F(R)t)^{2j+1-\frac{4k}{n}-2l_1} (1+G(R)t)^{1-2j-2l_2}t^{2l_1+2l_2}\\
	\times \prod_{i_1=1}^{l_1}F^{(1+m_{1i_1})}(R)^2\prod_{i_2=1}^{l_2}G^{(1+m_{2i_2})}(R)^2dR <\I
\end{multline}
for each $k \in [2,\ceil{s}+4]$, $j \in [0,k]$,
$l_i\ge0$ ($i=1,2,3$) with $l_2<j$ and $l_1+l_2+l_3\le j$, 
and $m_i \in (\N\cup\{0\})^{l_i}$ ($i=1,2$) with $|m_1|+|m_2|=j-l_1-l_2-l_3$.
We divide $\int_0^\I=\int_0^{R_0}+\int_{R_0}^{R_1}+\int_{R_1}^\I$
and denote the left hand side of \eqref{eq:pf_regv1} as $J_1+J_2+J_3$.
As in Step 2-b,
it is easy to see that $J_2$ is finite and tends to zero as $t\to\I$.

We now estimate $J_1$.
By \eqref{eq:decayav0}, \eqref{eq:decayFG0}, and \eqref{eq:FGratio},
$J_1$ is bounded by
\begin{align*}
	&C\int_0^{R_0} (R^{\kappa+1-l_3})^2 R^{-2k+2j+n-1} \prod_{i_1=1}^{l_1}(R^{\kappa-1-m_{1i_1}})^2
	\prod_{i_2=1}^{l_2}(R^{\kappa-1-m_{2i_2}})^2 dR \\
	={}&C\int_0^{R_0} R^{2\kappa-2(k-1)+n+2\kappa(l_1+l_2)-1} dR<\I
\end{align*}
because $2\kappa-2(k-1)+n > 2\kappa-2(\ceil{s}+3)+n>0$ by the choice of $\kappa$.

We finally treat $J_3$.
From \eqref{eq:decayFGinf}, it is bounded by
\begin{multline}\label{eq:estJ_3}
	C\int_{R_1}^\I \frac{R^{-n+2-2l_3-2\delta}}{ R^{2k-2j-n+1} }
	\prod_{i_1=1}^{l_1}(R^{-\frac{n}2-1-m_{1i_1}})^2\prod_{i_2=1}^{l_2}(R^{-\frac{n}2-1-m_{2i_2}-2\delta})^2\\
	\times \(\frac{1+F(R)t}{1+G(R)t}\)^{2j}(1+F(R)t)^{1-\frac{4k}{n}-2l_1} (1+G(R)t)^{1-2l_2}t^{2(l_1+l_2)}dR.
\end{multline}
By Young's inequality and \eqref{eq:decayFGinf}, we have
\[
	(1+F(R)t)^{1-\frac{4k}{n}-2l_1} (1+G(R)t)^{1-2l_2}t^{2(l_1+l_2)}
	\le C R^{n(l_1+l_2)+2\delta \max(0,2l_2-\frac{4k}{n}+1)}.
\]
Therefore, \eqref{eq:estJ_3} is bounded by
\[
	C\int_{R_1}^\I 
	R^{-2(k-1)+\delta(4j-2-2\min(2 l_2,\frac{4k}{n}-1))-1} dR,
\]
which is finite if $-2(k-1)+\delta(4j-2-2\min(2l_2,\frac{4k}{n}-1))<0$.
Since $j\in[0,k]$ and $l_2\in[0,j]$,
this condition is satisfied for all $j$ and $l_2$ if 
$-2(k-1)+\delta(4k-2-2\min(0,\frac{4k}{n}-1))<0$.
The latter condition can be written as
\[
	\delta < \min \( \frac12 - \frac1{4k-2}, \frac{n}{2(n-2)}-\frac{n}{2k(n-2)}\)
	\le \min \( \frac13, \frac14 + \frac{1}{2(n-2)}\) 
\]
thanks to the fact that the range of $k$ is $[2,\ceil{s}+4]$.
\smallbreak
The estimate of ${\bf a}$ is similar and so we omit the detail.
We only note that the $L^2$ norm is conserved:
By the explicit representation of $a_0$,
\begin{align*}
	\norm{{\bf a}(t)}_{L^2(\R^n)}
	&{}=\int_0^\I |a_0(t,X(t,R))|^2X(t,R)^{n-1}\d_R X(t,R)dR\\
	&{}=\int_0^\I |A_0(R)|^2R^{n-1}dR.
\end{align*}
\end{proof}
\begin{remark}
The similar proof shows $\nabla \Phi \in Y^{s+4}_{p,q}$ also for
$p>2^*$ and $q\in (2,2+4\delta/n]$ at the sacrifice of the uniform bound in time.
We only have to replace the bounds of $1+F(R)t$ and $1+G(R)t$ 
with rougher ones:
\[
	1 \le 1+F(R)t,\, 1+G(R)t \le 1+ t \sup_R (|F(R)|+|G(R)|).
\]
We note that also the assumption $A_0(r)\neq 0$ ($r>0$)
is needed only for the uniform boundedness in time.
\end{remark}

\section{Preliminaries for the proof of Theorem \ref{thm:WKB}}\label{sec:pre1}
\subsection{Some estimates}
We first recall a consequence of the Hardy-Littlewood-Sobolev inequality,
which can be found in \cite[Th.~4.5.9]{Hormander1} or
\cite[Lemma~7]{PGEP}: 
\begin{lemma}\label{lem:Zhidkov1}
If $\varphi \in \mathcal{D}^\prime(\R^n)$ is such that $\d_j \varphi
\in L^p (\R^n)$ ($j=1,\dots,n$) for some $p\in]1,n[$, 
then there exists a constant $c$ such that $\varphi-c \in L^q(\R^n)$,
with $1/p=1/q+1/n$. 
\end{lemma}
The next two lemmas can be found in \cite{BGDD-IUMJ,LaJFA}: 
\begin{lemma}[Commutator estimate]\label{lem:com}
Let $s \ge 0$ be a real number
and $k\ge 0$ be an integer such that $k \le s$.
There exists $C>0$ such that
\[
        \Lebn{\Lambda^s(fg) - f \Lambda^s g}{2} \le
        C(\Lebn{\nabla f}{\I} \Sobn{g}{s-1} +
        \|\nabla^k f\|_{H^{s-k}}\Lebn{g}{\I}). 
\]
\end{lemma}
\begin{lemma}
Let $s \ge 0$ be a real number
and $k\ge 0$ be an integer such that $k \le s$. There exists $C>0$ such that
\begin{equation}\label{eq:KP1}
	\Lebn{\Lambda^s(fg)}{2} \le C(\Sobn{f}{s}
	\Lebn{g}{\I} + \Lebn{f}{\I}\|\nabla^k g\|_{H^{s-k}}),
\end{equation}
for all $f \in H^s \cap L^\I$ and $g\in \dot{H}^k \cap\dot{H}^{s} \cap L^\I$,
and that
\begin{equation}\label{eq:KP2}
        \Lebn{\Lambda^s\nabla (fg)}{2} \le C(\Sobn{\nabla f}{s}
        \Lebn{g}{\I} + \Lebn{f}{\I}\Sobn{\nabla g}{s}),
\end{equation}
for all $  f,g \in \dot{H}^1 \cap\dot{H}^{s} \cap L^\I$.
\end{lemma}
The following lemma is a modification of the results in \cite{AC-SP,CM-AA}:
\begin{lemma}\label{lem:Hartree}
Let $n \ge 3$, $k\in \R_+$, and $s \in \R$.
Let $\gamma\in(0,n)$ satisfy $\frac{n}{2} -k < \gamma \le n-k$.
Then, there exists $C_s$ such that 
\[
        \norm{\lvert\nabla\rvert^k(|x|^{-\gamma} * f)}_{H^{s}} \le C_s
        (\norm{f}_{H^{s}} + \norm{f}_{L^{1}}),\quad \forall f \in L^{1}
        \cap H^{s}. 
\]
\end{lemma}
\begin{proof}
Since $\F |x|^{-\gamma} = C|\xi|^{-n+\gamma}$ for $\g\in(0,n)$, it holds that
\[
        \norm{|\nabla|^k(|x|^{-\gamma}*f)}_{H^{s}} 
        = C\norm{\Jbr{\xi}^{s} |\xi|^{-n+\gamma+k} \F f}_{L^2}.
\]
The high frequency part $(|\xi|>1)$ is bounded by $C\norm{f}_{H^{s}}$
if $-n+\gamma+k \le 0$.
On the other hand, the low frequency part $(|\xi| \le 1)$ is bounded by
\begin{align*}
        C\Lebn{\F f}{\I}\(\int_{|\xi| \le 1}|\xi|^{2(-n+\gamma+k)}d\xi\)^{\frac12}
        \le C\Lebn{f}{1} 
\end{align*}
if $2(-n+\gamma+k)>-n$, that is, if $\gamma>n/2-k$.
\end{proof}

\subsection{Local existence of the unique solution of \eqref{eq:QSP}}
We now give an existence result of the unique solution of \eqref{eq:QSP}.
With slight generalization of the nonlinearity,
let us consider the following system of Hartree type:
\begin{equation}\label{eq:QHE}
\left\{
	\begin{aligned}
		&\partial_t a^\eps + \nabla \phi^\eps \cdot \nabla a^\eps + \frac12 a^\eps \Delta \phi^\eps =i\frac{\eps}{2} \Delta a^\eps, 
		&&a^\eps(0,x)=A_0^\eps(x); \\
		&\partial_t \phi^\eps +\frac12 |\nabla \phi^\eps|^2 + \l (|x|^{-\g}*|a^\eps|^2) = 0 ,
		&&\phi^\eps(0,x)= \Phi_0(x).
	\end{aligned}
\right.
\end{equation}
The system \eqref{eq:QSP} corresponds to 
\eqref{eq:QHE} with $\g=n-2$. 
\begin{assumption}\label{asmp:existenceHE}
Let $n \ge 3$ and $\max(n/2-2,0) < \gamma \le n-2$. Let $\l \in \R$.
We suppose the following conditions with some $s>n/2+1$:\\
$\bullet$ The initial amplitude $A_0^\varepsilon \in H^{s+1}(\R^n)$ uniformly in
 $\varepsilon \in (0,1]$. \\
$\bullet$  The initial phase $\Phi_0 \in C^4(\R^n)$ satisfies
$\nabla \Phi_0 \in Y^{s+2}_{p,q}(\R^n)$ for some $p\in (2^*,\I]$ and $q\in (2,n)$ with $p\ge q$.
\end{assumption}
\begin{theorem}\label{thm:existenceHE}
Let Assumption~\ref{asmp:existenceHE} be satisfied.
Then, there exists $T>0$ independent of $\varepsilon$ and $s$
such that, for all $\eps\in (0,1]$, \eqref{eq:QHE} has a unique solution 
\begin{equation*}
  (a^\varepsilon,\phi^\varepsilon)\in C([0,T]; H^{s+1}(\R^n)\times C^4(\R^n))
\end{equation*}
with $\nabla \phi^\eps \in C([0,T];Y^{s+2}_{p,q})$.
Moreover, $u^\eps=a^\eps e^{i{\phi^\eps}/{\eps}}$ solves
\eqref{eq:HE} and
the solution $ (a^\varepsilon,\nabla \phi^\varepsilon)$ is bounded
in $L^\I([0,T]; H^{s+1}(\R^n)\times Y^{s+2}_{p,q}(\R^n))$ uniformly in $\eps \in (0,1]$ and the following properties hold:
\begin{itemize}
\item $\phi^\eps-\Phi_0 \in L^{\max(\frac{p}2,\frac{n}\g+)}(\R^n) \cap L^\I(\R^n)$.
\item If $n\ge 5$ and if $p$ in Assumption \ref{asmp:existenceHE} satisfies $2^*<p<n$, then $\Phi_0, \phi^\eps \in L^\I(\R^n)$.
\end{itemize}
\end{theorem}
\begin{remark}
$2^*=2n/(n-2)<n$ if and only if $n\ge5$.
\end{remark}
Denoting $v^\eps:=\nabla \phi^\eps$,
we obtain the following system:
\begin{equation}\label{eq:system1}
\left\{
  \begin{aligned}
    &\partial_t a^\eps + v^\eps \cdot \nabla
	a^\eps + \dfrac{1}{2}a^\eps \nabla\cdot v^\eps 
	= i \frac{\eps}{2} \Delta a^\eps,&&
	a^\eps_{|t=0} = A_0^\eps, \\ 
	&\d_t v^\eps + v^\eps \cdot \nabla
	v^\eps + \l \nabla \(|x|^{-\g}
	\ast |a^\eps|^2\) = 0, && v^\eps_{|t=0}=\nabla \Phi_0.
  \end{aligned}
\right.
\end{equation}
We first solve this system.
Then, as we seen below, we can reconstruct $\phi^\eps$ from $v^\eps$.
The proof goes along the classical energy argument.
Then, the main part of the proof is to establish an a priori estimate.
We hence perform precisely only this part.
As a first step, we show the following proposition.
\begin{proposition}\label{prop:partialenergyest}
Let Assumption \ref{asmp:existenceHE} be satisfied.
If $(a^\eps,v^\eps)$ is a solution of \eqref{eq:system1} in
$C([0,T];H^{s+1}\times Y^{s+2}_{p,q})$,
then its ``partial energy'' $E_{\mathrm{part}}(t):=\Sobn{a^\eps}{s+1}^2 + \Sobn{\nabla^2 v^\eps}{s}^2$
satisfies
\[
	\frac{d}{dt}E_{\mathrm{part}}(t) \le C E_{\mathrm{part}}(t)^{\frac32}.
\] 
\end{proposition}
\begin{proof}
We first estimate the $H^{s+1}$ norm of $a^\eps$.
We use the following convention for the scalar product in $L^2$:
\[
        \Jbr{\varphi, \psi} := \int_{\R^n} \varphi(x) \overline{\psi(x)} dx.
\]
The notation $\Lambda=(1-\Delta)^{1/2}$ is also used.
Then,
\[
        \frac{d}{dt} \Sobn{a^\eps}{s+1}^2 = 2\Re \Jbr{\partial_t
        \Lambda^{s+1} a^\eps, \Lambda^{s+1} a^\eps}. 
\]
Let us bound the right hand side with the relation
\begin{equation*}
	\partial_t \Lambda^{s+1} a^\eps + \Lambda^{s+1} ( v^\eps \cdot \nabla a^\eps) 
	+ \frac{1}{2} \Lambda^{s+1}  ( a^\eps \nabla \cdot v^\eps)
	-i \frac{\eps}{2} \Delta \Lambda^{s+1}
	a^\eps= 0. 
\end{equation*}
This part is standard (for details, see \cite{AC-SP,CM-AA}). 
The point is that we cannot not use $\Lebn{\nabla v^\eps}{2}$ as a bound.
This is done by the use of Lemma \ref{lem:com} and
\eqref{eq:KP1} with suitable $k$.
For example, Lemma \ref{lem:com}
with $k=2$ shows the estimate
\begin{multline*}
	|\Re \Jbr{[\Lambda^{s+1}, v^\eps]\cdot  \nabla a^\eps,
	\Lambda^{s+1} a^\eps }|\\ 
	\le C(\Lebn{ \nabla v^\eps}{\I} \Sobn{\nabla a^\eps}{s}
	+ \Sobn{ \nabla^2 v^\eps}{s-1} \Lebn{\nabla a^\eps}{\I})
	\Sobn{a^\eps}{s+1},
\end{multline*}
in which $\Lebn{\nabla v^\eps}{2}$ does not appear.
As a result, we obtain
\begin{equation*}
	\frac{d}{dt} \Sobn{a^\eps}{s+1}^2
	\le C (\norm{a^\eps}_{W^{1,\infty}} + \norm{\nabla v^\eps}_{L^{\infty}})
	( \Sobn{a^\eps}{s+1} + \Sobn{\nabla^2 v^\eps}{s})\Sobn{a^\eps}{s+1}.
\end{equation*}
Recall that $ v^\eps \in Y^{s+2}_{p,q}$ 
and so that $\nabla v^\eps \to 0$ as $|x|\to\I$ by the definition of $Y$.
Hence, by the Sobolev embedding,
$\Lebn{\nabla v^\eps}{\I} \le C\Sobn{\nabla^2 v^\eps}{s}$.
We end up with
\begin{equation}\label{eq:unifbd6}
	\frac{d}{dt} \Sobn{a^\eps}{s+1}^2 \le C E_{\mathrm{part}}(t)^{\frac32}. 
\end{equation}
\smallbreak
Let us proceed to the estimate of $v^\eps$.
We denote the operator
$\Lambda^s \nabla^2$ by $Q$. 
We deduce from the equation for $v^\eps$ that
\begin{equation}\label{eq:unifbd7}
	\partial_t Q v^\eps +  Q ( v^\eps \cdot \nabla
	v^\eps) + Q \nabla (|x|^{-\g}*| a^\eps|^2)=0 .
\end{equation}
We consider the coupling of this equation and $Q v^\eps$.
The coupling with the second term of the left hand side of \eqref{eq:unifbd7}
can be written as
\begin{align*}
        \Jbr{Q ( v^\eps \cdot \nabla v^\eps), Q v^\eps}
        =& \Jbr{ v^\eps \cdot \nabla Q v^\eps, Q v^\eps} 
        + \Jbr{[\Lambda^s\nabla, v^\eps ] \cdot \nabla^2 v^\eps, 
		Qv^\eps} \\
        &+ \Jbr{\Lambda^s \nabla( \nabla v^\eps\cdot \nabla v^\eps), Qv^\eps}. 
\end{align*}
As the previous case,  integration by parts shows
\begin{equation}\label{eq:l:unifbf8}
        |\Re \Jbr{ v^\eps \cdot \nabla Q v^\eps, Q v^\eps}|
        \le \frac{1}{2} \Lebn{ \nabla v^\eps}{\I}
		\Sobn{ \nabla^2v^\eps }{s}, 
\end{equation}
and the commutator estimate with $k=1$ also shows
\begin{multline}\label{eq:unifbd9}
        |\Re \Jbr{[\Lambda^s\nabla , v^\eps]\cdot \nabla^2 v^\eps, 
		Qv^\eps}| \\
        \le C (\Lebn{ \nabla v^\eps}{\I} \Sobn{ \nabla^2 v^\eps}{s}
		+ \Sobn{\nabla^2 v^\eps}{s-1}\Lebn{ \nabla^2 v^\eps}{\I}
		)\Sobn{ \nabla^2 v^\eps}{s} .
\end{multline}
We also have
\begin{equation}\label{eq:unifbd95}
        |\Re \Jbr{\Lambda^s \nabla( \nabla v^\eps\cdot \nabla v^\eps), Qv^\eps}| 
        \le C \Lebn{ \nabla v^\eps}{\I} \Sobn{ \nabla^2 v^\eps}{s}
		\Sobn{ \nabla^2 v^\eps}{s} 
\end{equation}
by \eqref{eq:KP2}.
For the estimate of the Hartree nonlinearity, we use Lemma
\ref{lem:Hartree} with $k=2$ to obtain
\begin{align}
	\Sobn{\lambda \nabla^3 (|x|^{-\g}*|a^\eps|^2)}{s} 
	&\le C \Sobn{\nabla^2 (|x|^{-\g}*|a^\eps|^2)}{s+1} \nonumber \\
	&\le C (\Lebn{a^\eps}{\I} \Sobn{ a^\eps }{s+1}
	+ \Lebn{a^\eps}{2}^2) \label{eq:unifbd10}
\end{align}
as long as $\g\in (n/2-2.n-2]$.
Sum up \eqref{eq:unifbd7}--\eqref{eq:unifbd10} to conclude that
\begin{equation}\label{eq:unifbd11}
	\frac{d}{dt} \Sobn{\nabla^2 v^\eps}{s}^2 \le C E_{\mathrm{part}}(t)^{\frac32},
\end{equation}
which completes the proof.
\end{proof}
We now prove Theorem \ref{thm:existenceHE}.
\begin{proof}[Proof of Theorem \ref{thm:existenceHE}]
We first obtain the solution $(a^\eps,v^\eps)$  of \eqref{eq:system1} by the
energy method
and then integrate $v^\eps$ to construct $\phi^\eps$.
\smallbreak
{\bf Step 1}.
We shall show the existence of the solution
$(a^\eps,v^\eps) \in C([0,T];H^{s+1}\times Y^{s+2}_{p,q})$
of \eqref{eq:system1} for small $T>0$.
Let us derive an a priori estimate of the energy
\begin{equation}\label{eq:fullenergy}
	E(t): = \Sobn{a^\eps(t)}{s+1}^2 + \norm{v^\eps(t)}_{Y^s_{p,q}}^2.
\end{equation}
By Proposition \ref{prop:partialenergyest} and
Gronwall's lemma, there exists $T$ such that
\begin{equation}\label{eq:partenergybound}
	\sup_{t\in[0,T]}E_{\mathrm{part}}(t)
	\le 2E_{\mathrm{part}}(0).
\end{equation}
Next we estimate $v^\eps$ and $\nabla v^\eps$.
By the second equation of \eqref{eq:system1}, we obtain
\[
	v^\eps(t) = \nabla \Phi_0 - \int_0^t \((v^\eps \cdot \nabla)v^\eps
	+ \l \nabla (|x|^{-\g}*|a^\eps|^2)\)ds.
\]
Therefore, we deduce by the H\"older inequality that
\begin{equation}\label{eq:Lpestv}
\begin{aligned}
	\norm{ v^\eps}_{L^\I([0,T];L^p)} \le{}&
	\Lebn{\nabla \Phi_0}{p} + 
	T \norm{ v^\eps }_{L^\I([0,T];L^p)} 
	\norm{\nabla v^\eps}_{L^\I([0,T]\times \R^n)}\\
	&{} + T|\l|\norm{\nabla (|x|^{-\g}*|a^\eps|^2)}_{L^\I([0,T];L^p)}
\end{aligned}
\end{equation}
and
\begin{equation}\label{eq:Lpestdv}
\begin{aligned}
	\norm{\nabla  v^\eps}_{L^\I([0,T];L^q)} \le{}&
	\Lebn{\nabla^2 \Phi_0}{q} + 
	T \norm{\nabla  v^\eps }_{L^\I([0,T];L^q)} 
	\norm{\nabla v^\eps}_{L^\I([0,T]\times \R^n)}\\
	&{}+ T\norm{v^\eps}_{L^\I([0,T];L^p)}
	\norm{\nabla^2  v^\eps }_{L^\I([0,T];L^{\frac{pq}{p-q}})} \\
	&{} + T|\l|\norm{\nabla^2 (|x|^{-\g}*|a^\eps|^2)}_{L^\I([0,T];L^q)}.
\end{aligned}
\end{equation}
Notice that $H^s \hookrightarrow L^2 \cap L^\I \hookrightarrow L^{\frac{pq}{p-q}}$ since 
${pq}/({p-q})\in (2,\I]$ holds by assumption $p \ge q>2$.
Moreover, we infer from Lemma \ref{lem:Hartree} that
\begin{align*}
	&\Lebn{\nabla(|x|^{-\g}*|a^\eps|^2)}{p} \le 
	C \Lebn{|\nabla|^{1+n(\frac12-\frac1p)}(|x|^{-\g}*|a^\eps|^2)}{2} \le
	C \Sobn{a^\eps}{s+1}^2,\\
	&\Lebn{\nabla^2(|x|^{-\g}*|a^\eps|^2)}{q} \le 
	C \Lebn{|\nabla|^{2+n(\frac12-\frac1q)}(|x|^{-\g}*|a^\eps|^2)}{2} \le
	C \Sobn{a^\eps}{s+1}^2,
\end{align*}
provided $n/p-1<\g \le n-2$ and $n/q-2 < \g \le n-2$, respectively.
By the assumptions $p>2^*$ and $q>2$, 
$\max({n}/p-1 , {n}/q-2) < n/2-2$.
Letting $T$ so small that $T(2E_{\mathrm{part}}(0))<1/3$ if necessary,
one sees from \eqref{eq:partenergybound} that
\begin{equation}\label{eq:partenergybound2}
	\norm{ v^\eps}_{L^\I([0,T];L^p)} +
	\norm{\nabla  v^\eps}_{L^\I([0,T];L^q)} \le 
	3\Lebn{\nabla \Phi_0}{p}+ 3\Lebn{\nabla^2 \Phi_0}{q}
	+C(E_{\mathrm{part}}(0)).
\end{equation}
Plugging \eqref{eq:partenergybound2} to \eqref{eq:partenergybound},
we obtain the desired energy estimate:
There exist $T$ and $C$ such that
$\sup_{t\in[0,T]}E(t) \le C(E(0))$.
Thus, we see from a standard argument that a solution $(a^\eps,v^\eps)$
of \eqref{eq:system1}  
exists in $C([0,T];H^{s+1}\times Y^{s+2}_{p,q})$.
\smallbreak
{\bf Step 2}.
We next investigate the decay property of $v^\eps$
and show the uniqueness of the solution of \eqref{eq:system1}.
Since $q<n$ by assumption, $v^\eps\in Y^{s+2}_{p,q}(\R^n) \hookrightarrow L^p(\R^n) \cap L^\I(\R^n)$.
By the H\"older inequality and the Hardy-Littlewood-Sobolev inequality,
we have
\begin{equation}\label{eq:vtail}
	v^\eps-\nabla \Phi_0=- \int_0^t \((v^\eps \cdot \nabla)v^\eps
	+ \l \nabla (|x|^{-\g}*|a^\eps|^2)\)ds
	\in L^{\max\(\frac{pq}{p+q},\frac{n}{\g+1}+\)} \cap L^\I.
\end{equation}
Notice that $pq/(p+q)<p$ for all $p>2^*$
and that $n/(\g+1)<\min (2^*,n)<p$ since $\g>\max(n/2-2,0)$. 
Therefore, $v^\eps-\nabla \Phi_0$ decreases at spatial infinity faster than
$v^\eps$ and $\nabla \Phi_0$ themselves.

Let us proceed to the uniqueness of \eqref{eq:system1}.
Let $(a^\eps_1,v^\eps_1)$ and $(a^\eps_2,v^\eps_2)$ be two
solutions of \eqref{eq:system1} in $C([0,T];H^{s+1}\times Y^{s+2}_{p,q})$
with $(a_i^\eps,v_i^\eps)(0)=(A_0^\eps,\nabla \Phi_0)$.
Put $d_a^\eps=a_1^\eps-a_2^\eps$ and $d_v^\eps=v_1^\eps-v_2^\eps$.
We remark that $d_a^\eps(0)\equiv0$ and $d_v^\eps(0)\equiv 0$.
Moreover, we see from the above estimate \eqref{eq:vtail} that
$d_v^\eps=(v_1^\eps-\nabla \Phi_0)-(v_2^\eps-\nabla \Phi_0)$
and so $d_v^\eps\to 0$ as $|x|\to\I$.
Now, we estimate 
\[
	E_d(t):= \Lebn{d_a^\eps(t)}{2}^2 + \Lebn{\nabla d_v^\eps(t)}{2}^2.
\]
It is important to note that $\nabla v_1^\eps$ and $\nabla v_2^\eps$
do not necessarily belong to $L^2$ by definition of $Y^s_{p,q}$.
Nevertheless, their difference $d_v^\eps$ may do so because it is 
identically zero and so belongs to $L^2$ at the initial time.
By an energy estimate, we have
\begin{equation}\label{eq:unique5}
	\frac{d}{dt}E_d(t) \le C(\|a^\eps_i\|_{H^{s+1}},\|v^\eps_i\|_{Y^{s+2}_{p,q}})
	E_d(t).
\end{equation}
Hence, we conclude from Gronwall's lemma that
\[
	E_d(t) \le C(\|a^\eps_i\|_{H^{s+1}},\|v^\eps_i\|_{Y^{s+2}_{p,q}}) E_d(0)=0
\]
as long as the solutions $(a^\eps_i,v^\eps_i)$ exist.
This implies that $d_a^\eps\equiv0$ and $\nabla d_v^\eps\equiv0$.
In particular, there exists a function $d=d(t)$ of time such that
$d_v^\eps(t,x)=d(t)$.
Recall that $d_v^\eps(t,x) \to 0$ as $|x|\to\I$.
As a result, $d(t)\equiv0$ follows and
we hence obtain $(a_1^\eps,v^\eps_1)=(a_2^\eps,v^\eps_2)$.
\smallbreak
{\bf Step 3}.
We finally construct $\phi^\eps$ such that $\nabla \phi^\eps = v^\eps$.
Define $\phi^\eps$ by
\[
	\phi^\eps(t)=\Phi_0 - \int_0^t \( \frac12|v^\eps(s)|^2 + \l(|x|^{-\g}*|a^\eps|^2)(s) \)ds.
\]
Then, one easily verifies that $(a^\eps,\nabla \phi^\eps)$ solves \eqref{eq:system1}
and that $\nabla \phi^\eps \in Y^{s+2}_{p,q}$.
Since we have already known the uniqueness of the solution to \eqref{eq:system1}, 
$\nabla \phi^\eps=v^\eps$.
Thus, $(a^\eps,\phi^\eps)$ is a unique solution to \eqref{eq:QHE}.
Though $\phi^\eps$ and $\Phi_0$ themselves do not necessarily 
belong to any Lebesgue space, it follows from the 
H\"older inequality and the Hardy-Littlewood-Sobolev inequality that
\[
	\phi^\eps(t) -\Phi_0 = - \int_0^t \( \frac12|v^\eps|^2 + \l(|x|^{-\g}*|a^\eps|^2) \)ds \in L^{\max\(\frac{p}2,\frac{n}\g +\)}\cap L^\I.
\]
Moreover, it is bounded uniformly in $\eps\in[0,1]$.

If $n\ge 5$ and $2^*<p<n$ then, applying Lemma \ref{lem:Zhidkov1},
we see that there exists a constant $c_0 \in \R$ 
such that
$\|\Phi_0-c_0 \|_{L^{p^*}}\le C \|\nabla \Phi_0\|_{L^{p}}$.
Moreover, 
since $\Phi_0-c_0$ decays at the spatial inifinity,
it follows by the Sobolev inequality that $\Lebn{\Phi_0 - c_0}{\I} \le C \Sobn{\nabla^2 \Phi_0}{s}$,
which shows $\Phi_0 \in L^\I$ and so $\phi^\eps \in L^\I$.
Remark that $p/2<(p<)p^*$ and that $n/\g<n/(n/2-2)=2^{**}<p^*$
since $n\ge 5$ and $p>2^*$.
Therefore, the difference $\phi^\eps -\Phi_0$
decays faster than $\phi^\eps$ and $\Phi_0$.
\end{proof}

\begin{remark}
In \cite{AC-SP,CM-AA}, the key for existence result is to
solve the system for $(a^\eps , \nabla v^\eps)$ in $H^s
\times H^s$ first.
Here, we first solve the system for $(a^\eps , \nabla^2 v^\eps)$
in $H^{s+1}\times H^s$.
This is the difference.
The point is that even if $\nabla v^\eps \not\in L^2$,
we obtain the energy estimate by the Sobolev embedding:
If $n\ge3$ and $\nabla v^\eps \to 0$ as $|x|\to\I$ then
$\norm{\nabla v^\eps}_{L^\I} \le C \norm {\nabla^2 v^\eps}_{H^{s}}$.
We also note that it would be difficult to solve
the system for $(a^\eps,\nabla^3 v^\eps)$ in $H^{s+2}\times H^s$
if $n=3,4$ because $\nabla^3 v^\eps \in H^s$
do not yield this kind of bound on $\nabla v^\eps$, in general.
\end{remark}
\section{Proof of Theorem \ref{thm:WKB}}\label{sec:proof1}
We prove the following theorem.
\begin{theorem}\label{thm:globalexpansion}
Let Assumption \ref{asmp:existence} be satisfied.
Let $(a^\eps,\phi^\eps)$ be the solution to \eqref{eq:QSP}
given in Theorem \ref{thm:existenceHE}.
If \eqref{eq:limsys2} has a global solution $(a_0,\phi_0)$
which satisfies $\eta(T;s,p,q)<\I$ for all $T<\I$,
then there exist 
\[
	(a_j,\phi_j) \in C([0,\I); H^{s-2j+3}\times Y^{s-2j+5}_{2^*,2})
\]
($1\le j \le N$) and constant $C_s$ depending only on $n$ and $s$
such that, for any $T>0$, $(a^\eps,\phi^\eps)$ exists until $t=T$ and
it holds that
\begin{equation*}
\left\{
\begin{aligned}
	a^\eps&{}=a_0+\sum_{j=1}^N \eps^j a_j +O(\eps^{N+1}) \IN L^\I([0,T],H^{s-2N+1}(\R^n)), \\
	\phi^\eps&{}=\phi_0+\sum_{j=1}^N \eps^j \phi_j +O(\eps^{N+1}) \IN L^\I([0,T],Y^{s-2N+3}_{2^*,2}(\R^n))\\
\end{aligned}
\right.
\end{equation*}
for $\eps \le C \eta(T) e^{-3C_s\eta(T)T}$.
\end{theorem}
\begin{remark}
We note that $\phi^\eps$ 
itself does not necessarily belong to the space
$Y^{s+3}_{2^*,2}(\R^n)$
as shown in Theorem \ref{thm:existenceHE}.
\end{remark}
\begin{remark}
We will see from the following proof that
$\phi_j$ ($j\ge 1$) belongs to $C([0,\I);Y^{s-2j+5}_{[\frac{2p}{2+p},\I]\cap(1^{**},\I],2})$
and the above expansion of $\phi^\eps$ is valid in
$C([0,\I);Y^{s-2N+3}_{[\frac{2p}{2+p},\I]\cap(1^{**},\I],2})$.
Remark that if $n \ge 5$ then 
$Y^{s}_{[\frac{2p}{2+p},\I]\cap(1^{**},\I],2}\subset H^s$.
\end{remark}
Theorem \ref{thm:WKB} immediately follows from this theorem.
Notice that the main amplitude $\beta_0$ is not $a_0$ but $a_0 e^{i\phi_1}$.
There is an interaction between the amplitude part and the phase part
because of the presence of nonlinearity.
This fact also leads us to some ill-posedness results for the ``usual'', that is,
non-scaled nonlinear Schr\"odinger equations (\cite{AC-LR,CM-AA,ThJDE}).
Similarly, the function $\beta_j$ is defined by $(a_i,\phi_i)$ ($i=0,1,\cdots,j+1$).
We will see that $(a_j,\phi_j)$ ($j\ge1$) solves
the ``$j$-th linearized system'' of \eqref{eq:QSP}:
\begin{equation}\label{eq:QSPtembwLn1}
	\left\{
	\begin{aligned}
		&\d_t a_j + \sum_{i_1+i_2=j} \nabla\phi_{i_1} \cdot \nabla a_{i_2}  +
		\sum_{i_1+i_2=j}\frac{1}{2}  a_{i_1} \Delta \phi_{i_2} 
		- i\frac{1}{2}\Delta a_{j-1} =0, \\
		&\d_t \phi_j + \sum_{i_1+i_2=j} \frac12 \nabla\phi_i \cdot \nabla\phi_j +
		\l \sum_{i_1+i_2=j} (|x|^{-(n-2)}\ast \Re(a_{i_1}
		\overline{a_{i_2}})) = 0 . \\
		&a_{j}(0) = A_j, \quad \phi_{j}(0) = 0.
	\end{aligned}
	\right.
\end{equation}
We separate
the proof of Theorem \ref{thm:globalexpansion} into three steps,
\begin{description}
\item[zeroth order] estimate on $a^\eps-a_0$ and $\phi^\eps-\phi_0$,
\item[first order] estimate on $a^\eps-a_0 -\eps a_1$ and $\phi^\eps-\phi_0 -\eps \phi_1$,
\item[higher order] estimate on $a^\eps - \sum_{j=0}^k \eps^k a_k$ and $\phi^\eps - \sum_{j=0}^k \eps^k \phi_k$ for $k\ge2$.
\end{description}
However, we only prove the third step because,
in order to exclude the dependency of $C_s$ on the expansion level $N$,
the main step is the third part.
This constant $C_s$ is chosen later (in the proof of Proposition \ref{prop:higherb}, below).

\subsection{Proof the theorem -- part 1: the zeroth order}\label{subsec:proof1}
We first state the estimate on the differences $a^\eps-a_0$ and
$\phi^\eps-\phi_0$.
\begin{proposition}\label{prop:zeroth}
Let Assumption \ref{asmp:existence} be satisfied.
Let $(a^\eps,\phi^\eps)$ be the solution to \eqref{eq:QSP}
given by Theorem \ref{thm:existenceHE}
and $(a_0,\phi_0)$ be the global solution to \eqref{eq:limsys2}
with $\eta(T)<\I$ for all $T<\I$.
Then, there exists a constant $C_s$ depending on $n$ and $s$, 
and $\Gamma_1$ depending on $A_0^\eps$ such that,
for any large $T$,
\begin{equation}\label{eq:gerror}
	\norm{a^\eps-a_0}_{L^\I([0,T],H^{s+1})} + \norm{\nabla \phi^\eps - \nabla \phi_0}_{L^\I([0,T],H^{s+2})}
	\le \eps \Gamma_1 e^{C_s \eta (T) T}
\end{equation}
holds for all $\eps \le \eps_0(T) \le \eta(T)C e^{-C_s\eta(T)T}$.
In particular, the existence time $T$ of $(a^\eps,\phi^\eps)$
can be chosen so that $\eps \sim \eta(T)e^{-C_s\eta(T)T}$.
\end{proposition}
The proof proceeds as in \cite{LT-MAA}. 
\subsection{Proof of the theorem -- part 2: the first order}
We next claim the following two points:
First is that $(a_1,\phi_1)$ is defined globally in time as
a limit $\eps\to0$ of $(\tildeE{a}_0,\tildeE{\phi}_0)$ (Proposition
\ref{prop:defaphi1}).
Second is that the asymptotics
\[
	a^\eps=a_0 +\eps a_1 + O(\eps^2), \quad
	v^\eps=v_0 +\eps v_1 + O(\eps^2)
\]
holds for large time (Proposition \ref{prop:first}).
\begin{proposition}\label{prop:defaphi1}
Let Assumption \ref{asmp:existence} be satisfied.
Suppose that \eqref{eq:limsys2} has a global solution $(a_0,\phi_0)$
which satisfies $\eta(T)<\I$ for all $T<\I$.
Then, there exists
$(a_1,\phi_1)\in C([0,\I),H^{s+1}\times
	Y^{s+3}_{2^*,2})$
which solves \eqref{eq:QSPtembwLn1}.
Let $E_1(t):=\Sobn{a_1(t)}{s+1}+\norm{\nabla \phi_1(t)}_{H^{s+2}}$.
Then, for any $T>0$, we have the following bound
\begin{equation}\label{eq:E1bound}
	\sup_{t\in[0,T]} E_1(t) \le \Gamma_1 e^{C_s\eta(T)T}
	=: \eta_1(T),
\end{equation}
where $\Gamma_1$, $C_s$, and $\eta$ are the same
as in Proposition \ref{prop:zeroth}.
\end{proposition}
\begin{proposition}\label{prop:first}
Let Assumption \ref{asmp:existence} be satisfied.
Let $(a^\eps,\phi^\eps)$ be the solution to \eqref{eq:QSP}
given by Theorem \ref{thm:existenceHE}.
Suppose that \eqref{eq:limsys2} has a global solution $(a_0,\phi_0)$
which satisfies $\eta(T)<\I$ for all $T<\I$.
Let $(a_1,\phi_1)$ be the limit defined in Proposition \ref{prop:defaphi1}.
Let $C_s$ be the same one as in Proposition \ref{prop:zeroth}.
Then, there exists a constant 
$\Gamma_2$ depending on $A_0^\eps$ such that
\begin{equation}\label{eq:gerror1}
\begin{aligned}
	\norm{a^\eps-a_0- \eps a_1}_{L^\I([0,T],H^{s-1})} + \norm{\nabla( \phi^\eps - \phi_0- \eps\phi_1)}_{L^\I([0,T],H^{s})} \\
	\le \eps^2 \Gamma_2 \eta(T)^{-1} e^{3C_s \eta (T) T}
	e^{\eps C_s \eta_1 (T) T}
\end{aligned}
\end{equation}
holds for all $0<\eps \le \eps_1(T) \le C e^{-2C_s\eta(T)T}$.
In particular, the existence time $T$ of $(a^\eps,\phi^\eps)$
can be chosen so that $\eps \sim e^{-2C_s\eta(T)T}$.
\end{proposition}
They are shown as in Propositions \ref{prop:highera} and \ref{prop:higherb}, below, respectively.
\subsection{Proof of the theorem -- part 3: higher order}\label{subsec:proof3}
We now consider the higher order expansion.
Assume that the constant $N$ in Assumption \ref{asmp:existence}
is bigger than one.
It is because if $N=1$ then the proof of Theorem \ref{thm:globalexpansion} is already finished with
Proposition \ref{prop:first}.
The proof is based on the induction argument.
We make following notation and definitions:
Our goal is to show that the asymptotics
\begin{equation}\label{eq:destinationm}
\left\{
\begin{aligned}
	a^\eps&{}=a_0+\sum_{j=1}^m \eps^j a_j +O(\eps^{m+1}) \IN L^\I([0,T],H^{s-2m+1}(\R^n)), \\
	\phi^\eps&{}=\phi_0+\sum_{j=1}^m \eps^j \phi_j +O(\eps^{m+1}) \IN L^\I([0,T],Y^{s-2m+3}_{2^*,2}(\R^n))\\
\end{aligned}
\right.
\end{equation}
for $m=N$. We define the following function:
\begin{equation}\label{def:etak}
	\eta_j(T):=\frac{\Gamma_j}{\eta(T)^{j-1}} e^{(2j-1)C_s \eta(T) T}
\end{equation}
with $\eta(T)$ is an increasing function defined in \eqref{def:eta},
$\Gamma_1$ and $\Gamma_2$ are as in Propositions
\ref{prop:zeroth} and  \ref{prop:first}, respectively,
and $\Gamma_j$ ($j\ge3$) is a constant depending only on $A_0^\eps$
to be chosen later. Note that if $T\gg1$ then
\begin{equation*}
	\eta_m(T) \gg \eta_{m-1}(T) \gg \dots \gg \eta_1(T) \gg \eta(T)>0.
\end{equation*}
The following two propositions complete the proof of Theorem \ref{thm:globalexpansion}.
\begin{proposition}\label{prop:highera}
Let Assumption \ref{asmp:existence} be satisfied for some $N\ge2$.
Suppose that \eqref{eq:limsys2} has a global solution $(a_0,\phi_0)$
which satisfies $\eta(T;s,p,q)<\I$ for all $T<\I$.
Fix $k_0 \in [1, N-1]$.
Assume that 
\[
	(a_j,\phi_j) \in C([0,\I); H^{s-2j+3} \times
Y^{s-2j+5}_{2^*,2})
\]
 ($1\le j \le k_0$) exist and all of them solve \eqref{eq:QSPtembwLn1}.
We further assume that there exists a positive constant $\Gamma_{k_0+1}$ such that
\[
	\overline{\lim_{\eps\to0}}\sup_{t\in[0,T]}\(\Sobn{\frac{a^\eps -\sum_{j=0}^{k_0}\eps^j a_j}{\eps^{k_0+1}}}{s-2{k_0}+1} + \norm{
	\frac{\nabla(\phi^\eps -\sum_{j=0}^{k_0}\eps^j \phi_j)}{\eps^{k_0+1}}}_{H^{s-2k_0+2}}\)
\]
is bounded by $\eta_{k_0+1}(T)$ defined in \eqref{def:etak}
for any fixed $T>0$.
Then, there exists $(a_{k_0+1},\phi_{k_0+1})
\in C([0,\I); H^{s-2k_0+1} \times Y^{s-2k_0+3}_{2^*,2})$
which solves \eqref{eq:QSPtembwLn1} and satisfies
\[
	\sup_{t\in[0,T]}\(
	\Sobn{a_{k_0+1}}{s-2{k_0}+1} + \norm{
	\nabla \phi_{k_0+1}}_{H^{s-2k_0+2}}
	\)\le \eta_{k_0+1}(T).
\]
\end{proposition}
\begin{proposition}\label{prop:higherb}
Let Assumption \ref{asmp:existence} be satisfied for some $N\ge2$.
Suppose that \eqref{eq:limsys2} has a global solution $(a_0,\phi_0)$
which satisfies $\eta(T;s,p,q)<\I$ for all $T<\I$.
Fix $k_0 \in [1,N-1]$.
Assume that, for all $1\le j \le k_0+1$, the solution
$(a_j,\phi_j) \in C([0,\I); H^{s-2j+3} \times
Y^{s-2j+5}_{2^*,2})$
of  \eqref{eq:QSPtembwLn1} exists
and satisfies
\[
	\sup_{t\in[0,T]}
	\(\Sobn{a_{j}}{s-2j+3} + \norm{
	\nabla \phi_{j}}_{H^{s-2j+4}}\)
	\le \eta_{j}(T).
\]
Then,  for any fixed $T>0$,
\[
	\sup_{t\in[0,T]} \(\Sobn{\frac{a^\eps -\sum_{j=0}^{k_0+1}\eps^j a_j}{\eps^{k_0+2}}}{s-2{k_0}-1} + \norm{
	\frac{\nabla(\phi^\eps -\sum_{j=0}^{k_0+1}\eps^j \phi_j)}{\eps^{k_0+2}}}_{H^{s-2k_0}}\)
\]
is bounded uniformly in $\eps \in (0,\eps_{k_0+2}]$.
In particular, the asymptotics \eqref{eq:destinationm} holds
with $m=k_0+1$ for $\eps \in (0,\eps_{k_0+2}]$.
$\eps_{k_0+2}$ can be chosen so that
$\eps_{k_0+2}\le C \eta(T) e^{-3C_s\eta(T)T}$.
Moreover, there exists a positive constant $\Gamma_{k_0+2}$
depending only on $A_0^\eps$
such that $\eta_{k_0+2}(T)$ defined in \eqref{def:etak} bounds
\[
	\overline{\lim_{\eps\to0}}\sup_{t\in[0,T]} \(\Sobn{\frac{a^\eps -\sum_{j=0}^{k_0+1}\eps^j a_j}{\eps^{k_0+2}}}{s-2{k_0}-1} + \norm{
	\frac{\nabla(\phi^\eps -\sum_{j=0}^{k_0+1}\eps^j \phi_j)}{\eps^{k_0+2}}}_{H^{s-2k_0}}\)
\]
for any fixed large $T>0$.
\end{proposition}
Proposition \ref{prop:first} implies that the
assumption of Proposition \ref{prop:highera} is satisfied for $k_0=1$.
Then, we see by induction that Proposition \ref{prop:higherb} holds for $k_0=N-1$.
Then, this gives \eqref{eq:destinationm} with $m=N$.
Before the proof, we introduce some more notation.
We write
\begin{align*}
		b_{m}^\eps &= \frac{a^\varepsilon
        -\sum_{j=0}^{m}\varepsilon^j  a_j}{\varepsilon^{m+1}}, & 
		w_{m}^\eps &= \frac{\nabla \phi^\varepsilon
        -\sum_{j=0}^{m}\eps^j  \nabla  \phi_j}{\varepsilon^{m+1}}. 
\end{align*}
An elementary computation shows that
$(b_{m}^\eps,w_{m}^\eps)$ satisfies
\begin{multline}\label{eq:QSPtembwLm1}
	\partial_t b_{m}^\eps
		+ \eps^{m+1}\(w_m^\eps \cdot\nabla
		b_{m}^\eps + \frac{1}{2}b_{m}^\eps\nabla \cdot w_{m}^\eps \) \\
		+\sum_{\ell=0}^{m}\eps^\ell  \(  
		w_{m}^\eps\cdot\nabla a_\ell + v_\ell\cdot\nabla b_{m}^\eps 
		+ \frac{1}{2}b_{m}^\eps \nabla \cdot v_\ell
		+ \frac{1}{2}a_\ell \nabla \cdot v_{m}^\eps\right) \\ 
		+ \sum_{\ell=0}^{m-1}\eps^{\ell}
		\sum_{i,j\le m, i+j=m+1+\ell} \left(v_i \cdot \nabla
		a_j  + \frac{1}{2}  a_i \nabla \cdot v_j\right)                 
		- i\frac{1}{2}\Delta a_{m} =
		i\frac{\eps}{2}\Delta b_{m}^\eps,
\end{multline}
\begin{multline}\label{eq:QSPtembwLm2}
	\partial_t w^\eps_m 
	+ \eps^{m+1}\( w_{m}^\eps \cdot \nabla
	w_{m}^\eps +  \l \nabla(|x|^{-(n-2)}\ast |b_{m}^\eps|^2)\) \\
	+ \sum_{\ell=0}^{m}\eps^\ell (\left(
	w_{m}^\eps\cdot\nabla v_\ell + v_\ell\cdot\nabla
	w_{m}^\eps \right) + 
	\lambda \nabla(|x|^{-(n-2)}\ast \Re(a_\ell
	\overline{b_{m}^\eps})) \\ 
	+ \sum_{\ell=0}^{m-1}\eps^{\ell}
	\sum_{i,j\le m,i+j=m+1+\ell}  \left(v_i \cdot \nabla v_j  
	+ \lambda \nabla(|x|^{-(n-2)}\ast \Re(a_i
	\overline{a_j}))\right) = 0, 
\end{multline}
and
\begin{align}\label{eq:QSPtembwLm3}
        b^\eps_{m}(0) =& \sum_{j=0}^{k-1-m}
        \varepsilon^j A_{j+m+1} + \varepsilon^{k-m} r^\varepsilon_{k+1} , &
        w^\eps_{m}(0) &= 0
\end{align}
as long as $(a_0,v_0):=(a_0,\nabla \phi_0)$ and
$(a_j,v_j):=(a_j,\nabla \phi_j)$ ($1\le j \le m$) solve 
\eqref{eq:limsys2} and \eqref{eq:QSPtembwLn1}, respectively,
where $r_{k+1}^\eps$ is $\eps^{-{k+1}}(A_0^\eps-\sum_{j=0}^k
\eps^j A_j)$.
If Assumption \ref{asmp:existence} is satisfied then $r_{k+1}^\eps$ is
bounded in $H^{s+1}$ as $\eps\to0$.

\begin{proof}[Proof of Proposition \ref{prop:highera}]
By assumption, $(b_{k_0}^\eps,w_{k_0}^\eps)$ is uniformly bounded in
$L^\I([0,T), H^{s-2k_0+1}\times H^{s-2k_0+2})$
in the limit $\eps\to0$.
Therefore, extracting a subsequence, there exists a weak limit,
denoted by $(a_{k_0+1},v_{k_0+1})$, in the same class.
Moreover, we obtain the bound
\[
	\sup_{t\in[0,T]}\(
	\Sobn{a_{k_0+1}}{s-2{k_0}+1} + \norm{
	v_{k_0+1}}_{H^{s-2k_0+2}}
	\)\le \eta_{k_0+1}(T).
\]
by the lower semi-continuity of the weak limit.
Since $(b_{k_0}^\eps,w_{k_0}^\eps)$ solves
\eqref{eq:QSPtembwLm1}--\eqref{eq:QSPtembwLm3}, 
we see that $(a_{k_0+1},v_{k_0+1})$ solves
\begin{equation}\label{eq:QSPtemavLn1}
	\left\{
	\begin{aligned}
		&\d_t a_j + \sum_{i_1+i_2=j} v_{i_1} \cdot \nabla a_{i_2}  +
		\sum_{i_1+i_2=j}\frac{1}{2}  a_{i_1} \nabla\cdot v_{i_2} 
		- i\frac{1}{2}\Delta a_{j-1} =0, \\
		&\d_t v_j + \nabla\sum_{i_1+i_2=j} \frac12 v_i \cdot v_j +
		\l \nabla\sum_{i_1+i_2=j} (|x|^{-(n-2)}\ast \Re(a_{i_1}
		\overline{a_{i_2}})) = 0,  \\
		&a_{j}(0) = A_j, \quad v_{j}(0) = 0.
	\end{aligned}
	\right.
\end{equation}
for $j=k_0+1$.
By the way, once we know $(a_j,v_j)$ ($j=[0,k_0]$),
we can solve this system directly by a standard argument
and obtain unique solution $(a_{k_0+1},v_{k_0+1})$ in the same space.
Therefore, the above weak limit is the unique solution to \eqref{eq:QSPtemavLn1}.
We now define $\phi_{k_0+1}$ by
\[
	\phi_{k_0+1}(t)=-\int_0^t
	\(\sum_{i_1+i_2=j} \frac12 v_i \cdot v_j +
	\l \sum_{i_1+i_2=j} (|x|^{-(n-2)}\ast \Re(a_{i_1}
	\overline{a_{i_2}}))\) ds.
\]
Then, $\nabla \phi_{k_0+1}=v_{k_0+1}$ holds by the uniqueness
of \eqref{eq:QSPtemavLn1}.
Hence, $\nabla \phi_{k_0+1}$ is the unique solution to \eqref{eq:QSPtembwLn1} for $j=k_0+1$.
Since $v_0\in L^p$ and $v_j \in L^2$ ($j\ge1$),
\[
	\phi_{k_0+1} \in C([0,T];{Y^{s-2k_0+3}_{[\frac{2p}{2+p},\I]\cap(1^{**},\I],2}})	.
\]
$T$ is arbitrary, and so we obtain the proposition.
\end{proof}
\begin{proof}[Proof of Proposition \ref{prop:higherb}]
By assumption, we can define $(b_{k_0+1}^\eps,w^\eps_{k_0+1})$
solving \eqref{eq:QSPtembwLm1}--\eqref{eq:QSPtembwLm3}.
We will bound
\[
	\widetilde{E}_{k_0+1}(t) := \Sobn{b_{k_0+1}^\eps(t)}{s-2k_0-1}
	+ \norm{w_{k_0+1}^\eps(t)}_{ H^{s-2k_0}}.
\]
Recall that the quadratic part of
\eqref{eq:QSPtembwLm1}--\eqref{eq:QSPtembwLm2}
is the same of \eqref{eq:QSP} up to a constant,
and that the linear part of \eqref{eq:QSPtembwLm1}--\eqref{eq:QSPtembwLm2}
is essentially the same form. 
Hence, mimicking the estimates in the proof of Theorem \ref{thm:existenceHE},
we deduce that,
for any fixed $T>0$,
\begin{equation}\label{eq:LTestk2}
	\frac{d}{dt}\widetilde{E}_{k_0+1}(t)\le C_s(\eps^{k_0+1}
	\widetilde{E}_{k_0+1}(t)^2 + \mu_{k_0+1}^\eps \widetilde{E}_{k_0+1}(t)
	+ c_{k_0+1} \nu_{k_0+1}^\eps)	
\end{equation}
holds for all $t\in[0,T]$.
Here, we have used two functions: First is
\[
	\mu_{k_0+1}^\eps=\mu_{k_0+1}^\eps(T) :=  \eta(T) + \sum_{j=1}^{k_0+1} \eps^j \eta_j(T)
\]
which bounds the linear part
\[
	\sup_{t\in[0,T]}\(\Sobn{\sum_{\ell=0}^{k_0+1} \eps^\ell a_\ell}{s-2k_0+1}
	+ \norm{v_0}_{ Y^{s-2k_0+2}_{p,q}}
	+ \norm{\sum_{\ell=1}^{k_0+1} \eps^\ell v_\ell}_{ H^{s-2k_0}}\)
\]
and second is
\[
	\nu_{k_0+1}^\eps=\nu_{k_0+1}^\eps(T) := \eta_{k_0+1}(T) + \sum_{\ell=0}^{k_0} \eps^\ell \sum_{i=\ell+1}^{k_0+1} \eta_i(T)\eta_{k_0+2+\ell-i}(T)
\]
which is an upper bound of the constant terms
\begin{multline*}
	\sup_{t\in[0,T]} \Bigg( \frac{1}{2}\Sobn{\Delta a_{k_0}}{s-2k_0-1}
	+C \sum_{\ell=0}^{k_0}\eps^{\ell}
	\sum_{{1\le i,j\le k_0+1 ,}\atop{i+j=k_0+2+\ell}}  \\
	\(\norm{v_i}_{ H^{s-2k_0+2}}\norm{v_j}_{ H^{s-2k_0+2}}  
	+ \Sobn{a_i}{s-2k_0+1} \Sobn{a_j}{s-2k_0+1}\)\Bigg)
\end{multline*}
up to an adjusting constant $c_{k_0+1}$.
The constant $C_s$ comes from \eqref{eq:LTestk2}.
This constant is independent of $k_0$ because 
it has been already taken into account when we
use $\mu^\eps_{k_0+1}$ and $\nu^\eps_{k_0+1}$.
\smallbreak
We now show that $\sup_{t\in[0,T]}\widetilde{E}_{k_0+1}(t)$ is uniformly 
bounded for small $\eps$.
We keep fixing $T>0$.
By Assumption \eqref{asmp:existence},
we see that there exists a positive constant $\beta_{k_0+1}$ 
depending only on $A_0^\eps$ such that
$\widetilde{E}_{k_0+1}(0) \le \beta_{k_0+1}$
holds for $\eps \in (0,1]$.
Set a function
\[
	Z_{k_0+1}(t):=\widetilde{E}_{k_0+1}(t)\exp(-C_s \mu_{k_0+1}^\eps(T) t)
\]
and two constants
\begin{align*}
\delta_{k_0+1}&{}:=(1+\sqrt{1+\beta_{k_0+1}})^{-1},\\
	\theta_{m+1} &{}:= \frac{\delta \mu_{k_0+1}^\eps(T)}{2c_{k_0+1}\nu_{k_0+1}^\eps(T) (1-e^{-C_s \mu_{k_0+1}^\eps(T) T})}.
\end{align*}
Then, multiplying the both sides of \eqref{eq:LTestk2}
by $\frac{\theta_{k_0+1}\exp(-C_s t \mu_{k_0+1}^\eps )}{(1+\theta_{k_0+1} Z_{k_0+1}(t))^2}$, we obtain
\[
	\frac{\theta_{k_0+1} Z_{k_0+1}^\prime(t)}{(1+\theta_{k_0+1} Z_{k_0+1}(t))^2} \le
	C_s \eps^{k_0+2} e^{C_s t\mu_{k_0+1}^\eps} \theta_{k_0+1}^{-1}
	+ C_s c_{k_0+1} \nu_{k_0+1}^\eps e^{-C_s t\mu_{k_0+1}^\eps }\theta_{k_0+1},
\]
where we denote $\mu_{k_0+1}^\eps(T)$ and $\nu_{k_0+1}^\eps(T)$ by
$\mu_{k_0+1}^\eps$ and $\mu_{k_0+1}^\eps$,
respectively, for short.
Integration over $[0,t]$ gives
\begin{multline}\label{eq:Zmbound0}
	\frac{1}{1+\theta_{k_0+1} Z_{k_0+1}(t)} \ge \frac{1}{1+\theta_{k_0+1} \widetilde{E}_{k_0+1}(0)}\\
	-\frac{c_{k_0+1}\nu_{k_0+1}^\eps}{\mu_{k_0+1}^\eps} (1-e^{-C_st \mu_{k_0+1}^\eps})\theta_{k_0+1}
	-\frac{\eps^{k_0+2}}{\mu_{k_0+1}^\eps} (e^{C_s t\mu_{k_0+1}^\eps }-1) \theta_{k_0+1}^{-1}.
\end{multline}
Let us show that the right hand side of \eqref{eq:Zmbound0} is bounded by $\delta_{k_0+1}/2$ from below.
For simplicity, in the following, we omit the index $k_0+1$ and
denote $\beta_{k_0+1}$, $c_{k_0+1}$, $\delta_{k_0+1}$,
$\mu_{k_0+1}^\eps$, $\nu_{k_0+1}^\eps$, and
$\theta_{k_0+1}$ by $\beta$, $c$, $\delta$, 
$\mu^\eps$, $\nu^\eps$, and $\theta$, respectively.
We also omit $T$ variable in $\eta(T)$ and $\eta_j(T)$.
By the fact that $\eta_{j+1} \gg \eta_j$ for each $j$ and large $T$
and by definitions of $\mu^\eps $ and $\nu^\eps$,
$c\nu^\eps \ge \mu^\eps$ holds for all $\eps\in[0,1]$ if $T$ is large.
Then, replacing $T$ with larger one if necessary, we obtain
\begin{multline}\label{eq:Zmbound1}
	\frac{1}{1+\theta \widetilde{E}_{k_0+1}(0)} -\delta \ge \frac{1}{1+\theta \beta} -\delta
	=\frac{e^{C_s \mu^\eps T}-1}{e^{C_s \mu^\eps T}(1+\frac{\mu^\eps}{c\nu^\eps}\frac{\delta\beta}{2})-1} - \delta \\
	\ge 
	\frac{e^{C_s \mu^\eps T}(2-2\delta -\delta^2\beta )-2+2\delta}{e^{C_s \mu^\eps T}(2+\delta\beta)-2} 
	\ge \frac{2-2\delta -\delta^2\beta }{2(2+\delta\beta)}
	= \frac{\delta}2,
\end{multline}
where we have used the relation $1-2\delta -\delta^2 \beta=0$.
Moreover, $\theta$ is the minimizer of the quantity
\[
	\frac{c\nu^\eps}{\mu^\eps} (1-e^{-C_s \mu^\eps T})\theta^2
	-\delta \theta
	+\frac{\eps^{k_0+2}}{\mu^\eps} (e^{C_s \mu^\eps T}-1)
\]
and so this quantity becomes less than or equal to zero if
\begin{equation}\label{eq:epsm1}
	\eps \le \(\frac{\delta^2 (\mu^\eps)^2 e^{C_s \mu^\eps T} }{
	c\nu^\eps (e^{C_s \mu^\eps T}-1)^2}\)^{\frac1{k_0+2}}.
\end{equation}
We now replace this condition with stronger but clearer one.
We first let $\eps$ be so small that
\begin{equation}\label{eq:epsm2}
	\eps \le \min_{j\in[1, k_0+1]} \(\frac{\eta}{\eta_j}\)^{\frac1j}=
	 \min_{j\in[1, k_0+1]} \frac{\eta}{\Gamma_j^{1/j} e^{(2-1/j)C_s \eta(T)T}}.
\end{equation}
For such $\eps$, we have $\mu^\eps \le (k_0+2)\eta$ and,
by definition of $\eta_j$ \eqref{def:etak}, 
\begin{align*}
	\nu^\eps ={}& \eta_{k_0+1} + 
	\sum_{\ell=0}^{k_0}\eps^{\ell}
	\sum_{i=\ell+1}^{k_0+1} \eta_i \eta_{ k_0+2+\ell-i} \\
	\le{}& \eta_{k_0+1}+ \widetilde{\Gamma}_1 \frac{e^{(2k_0+2)C_s \eta T}}{\eta^{k_0}} +
	\frac{e^{(2k_0+3)C_s\eta T}}{\eta^{k_0+1}}\sum_{\ell=1}^{k_0}
	\widetilde{\Gamma}_2\eta
	\( \frac{\eps }{\eta}e^{(2-1/\ell ) C_s \eta T} \)^{\ell} \\
	\le{}& \eta_{k_0+1}+ \widetilde{\Gamma}_1 \frac{e^{(2k_0+2)C_s \eta T}}{\eta^{k_0}} +
	\widetilde{\Gamma}_3 \frac{e^{(2k_0+3)C_s\eta T}}{\eta^{k_0}} 
	\le \widetilde{\Gamma}_4 \frac{e^{(2k_0+3)C_s \eta T}}{\eta^{k_0}},
\end{align*}
where $\widetilde{\Gamma}_i$ is a constant depending on $k_0$
and $\Gamma_{j}$ ($1\le j \le k_0+1$).
Therefore, the right hand side of \eqref{eq:epsm1} is bounded  below by
\begin{align*}
	\(\frac{\delta^2 (\mu^\eps)^2 e^{C_s \mu^\eps T} }{
	c\nu^\eps (e^{C_s \mu^\eps T}-1)^2}\)^{\frac1{k_0+2}}\ge {}&
	\(\frac{\delta^2 \eta^2 }{
	c\nu^\eps e^{C_s \mu^\eps T}}\)^{\frac1{k_0+2}} \\
	\ge{}&
	\widetilde{\Gamma}_5
	\(\frac{\eta^2 }{(\eta^{-k_0}e^{(2k_0+3)C_s\eta T})
	 e^{(k_0+2)C_s \eta T}}\)^{\frac1{k_0+2}} \\
	=&{}\widetilde{\Gamma}_5  \frac{\eta}{e^{(3-\frac{1}{k_0+2})C_s\eta T}} 
	 \ge \widetilde{\Gamma}_5 \frac{\eta}{e^{3C_s \eta T}}
	 =:\eps^{k_0+2},
\end{align*}
where $\widetilde{\Gamma}_5$ depends on
 $\widetilde{\Gamma}_4$, $\beta$, and $c$.
Then, the condition $\eps \le \eps_{k_0+2}$ ensures
\eqref{eq:epsm1} and so
\begin{equation}\label{eq:Zmbound2}
		\delta -\frac{c \nu^\eps}{\mu^\eps} (1-e^{-C_s \mu^\eps T})\theta
	-\frac{\eps^{k_0+2}}{\mu^\eps} (e^{C_s T\mu^\eps}-1) \theta^{-1}\ge 0.
\end{equation}
Note that $\eps_{k_0+2}$ is smaller than the right hand
side of \eqref{eq:epsm2} and so that $\eps\le \eps_{k_0+2}$
is stronger than \eqref{eq:epsm2}.

Furthermore, plugging \eqref{eq:Zmbound1} and \eqref{eq:Zmbound2}
to \eqref{eq:Zmbound0}, we obtain
\begin{equation}\label{eq:goalest}
	\sup_{t\in[0,T]} \widetilde{E}_{k_0+1}(t) \le 3 \sqrt{1+\beta} \theta^{-1} e^{C_s \mu^\eps T}
	\le \frac{6c\sqrt{1+\beta}\nu^\eps}{\delta \mu^\eps} e^{C_s \mu^\eps T},
\end{equation}
which is the desired bound.
Indeed, the right hand side is bounded by
\[
	\frac{6c\sqrt{1+\beta}\widetilde{\Gamma}_4 }{\delta \eta^{k_0+1}} e^{C_s (3k_0+5)\eta T}
\]
as long as $\eps\le\eps_{k_0+2}$.
We finally confirm that the right hand side of \eqref{eq:goalest}
tends to $\eta_{k_0+2}(T)$ defined by \eqref{def:etak} with a suitable constant.
It holds that
\begin{align*}
	\lim_{\eps\to0}\mu^\eps
	&{}=\lim_{\eps\to0}\mu^\eps_{k_0+1}(T)= \eta(T), \\
	\lim_{\eps\to0}\nu^\eps
	&{}= \lim_{\eps\to0}\nu^\eps_{k_0+1}(T)
	= \eta_{k_0+1} + \sum_{i=1}^{k_0+1} \eta_i \eta_{k_0+2-i} 
	\le \frac{\hat\Gamma_{k_0+2}}{\eta(T)^{k_0}}e^{(2k_0+2)C_s\eta(T)T},
\end{align*}
where $\hat\Gamma_{k_0+2}$ depends on $k_0$
and $\Gamma_j$ ($1\le j \le k_0+1$).
Therefore, we end up with the estimate 
\begin{align*}
	\limsup_{\eps\to0}\sup_{t\in[0,T]} \widetilde{E}_{k_0+1}(t) 
	\le{}& \frac{6\sqrt{1+\beta}(\hat\Gamma_{k_0+2}\eta(T)^{-k_0}e^{(2k_0+2)C_sg(T)T})}{\delta \eta(T)} e^{C_s \eta(T)T} \\
	=:{}&\frac{\Gamma_{k_0+2}}{ \eta(T)^{k_0+1}} e^{(2k_0+3)C_s \eta(T)T}=\eta_{k_0+2}(T),
\end{align*}
which completes the proof.
\end{proof}
\appendix
\section{Global existence of solution to the Euler-Poisson equations}\label{sec:EP}
In this section we give the proof of Theorem \ref{thm:EP}.
Let us consider
\begin{equation}\label{eq:a:rEP}
\left\{
\begin{aligned}
	&\rho_t + r^{-(n-1)}\d_r (r^{n-1} \rho v) = 0, 
	&&\rho(0,r) = \rho_0(r);\\
	&v_t + v \d_r v + \l \d_r V_{\mathrm{P}} =0,
	&&v(0,r) = v_0(r);\\
	&-r^{-(n-1)}\d_r (r^{n-1} V_{\mathrm{P}}) = \rho,\, V_{\mathrm{P}} \in L^\I, && V_{\mathrm{P}}\to 0\text{ as }r\to \I
\end{aligned}
\right.
\end{equation}
Here, $n\ge1$ denotes space dimensions,
$r\ge0$ denotes the distance from the origin, and $\l$ is a given physical constant.
Unknowns are the mass density $\rho=\rho(t,r)\ge 0$ and
 the velocity field $v=v(t,r) \in \R$.
If $n=1$ or $2$, we change the condition for Poisson equation of \eqref{eq:a:rEP}
into $V_{\mathrm{P}} (t,0)=0$, $\d_rV_{\mathrm{P}} \in L^\I$,
and $\d_rV_{\mathrm{P}}\to 0$ as $r\to \I$.
In other words, we let $V_{\mathrm{P}}$ be as
\begin{equation}\label{eq:asmp:V_p}
	V_{\mathrm{P}}(t,r) = 
	\left\{
	\begin{aligned}
	&\int_0^r s^{-(n-1)} \int_0^s \rho(t,\sigma)\sigma^{n-1} d\sigma ds &&\text{ if } n=1,2, \\
	&\int_r^\I s^{-(n-1)} \int_0^s \rho(t,\sigma)\sigma^{n-1} d\sigma ds &&\text{ if } n \ge 3.
	\end{aligned}
	\right.
\end{equation}
This is well-defined because we restrict our attention to $\rho$ belonging to
$C([0,\I)) \cap L^1((0,\I),r^{n-1}dr)$.
One can verify that the condition $V_{\mathrm{P}}(t,r)\to0$ as $r\to\I$ is not
suitable for $n=1$ or $2$.
Remark that \eqref{eq:a:rEP} is a radial version of
the compressible Euler-Poisson equations
\begin{equation}\label{eq:a:EP}
\left\{
  \begin{aligned}
    &\partial_t \rho + \mathrm{div} (\rho v)=0,&&
	\rho(0,x) = \rho_0(x), \\ 
	&\partial_t v + (v\cdot \nabla) v + \l \nabla V_{\mathrm{P}}  = 0, && v(0,x)=v_0, \\
    &-\Delta V_{\mathrm{P}} = \rho, \quad V_{\mathrm{P}}\in L^\I(\R^n),
    &&V_{\mathrm{P}} \to 0 \text{ as }|x|\to\I.
  \end{aligned}
\right.
\end{equation}

\subsection{Reduction to an ODE for characteristic curves}
We follow the argument in \cite{ELT-IUMJ,MaRIMS}.
Define characteristic curve $X$ as a function $\R_+ \to \R_+$ with parameter $R\in\R_+$
which is defined by an ODE
\[
	\frac{d}{dt}X(t,R) = v(t,X(t,R)), \quad X(0,R)=R.
\]
Let $m(t,r):=\int_0^r \rho(t,s)s^{n-1}ds$.
In the followings, we denote $\d_t X$ as $X^\prime$ by the respect that $R$ is a parameter.
Then, \eqref{eq:a:rEP} is reduced to the following ODE for characteristic curve
\begin{equation}\label{eq:ODE1}
	X^{\prime\prime}(t,R) = -\frac{\l m_0(R)}{X(t,R)^{n-1}},\quad
	X^\prime (0,R)=v_0(R) , \quad X(0,R)=R,
\end{equation}
where $m_0(r):=\int_0^r \rho_0(s)s^{n-1}ds$.
This reduction is the same spirit as the use of the Lagrangian coordinate
(see \cite{NN-IUMJ,NNY-MMMAS}).
Put 
\[
	B(t,R) := \exp \( \int_0^t \partial_r u(\tau, X(\tau,R)) d\tau \).
\]
It holds that $B(t,R)=\partial_R X(t,R)$.
The solution to \eqref{eq:a:rEP} is given explicitly in terms of $X$ and $B$ as
\begin{align}\label{eq:Xtorho}
	\rho(t,X(t,R))&{}=\frac{R^{n-1}\rho_0(R)}{X^{n-1}B(t,R)}, &
	v(t,X(t,R))&{}=\frac{d}{dt}X(t,R). 
\end{align}
As in \cite{MaRIMS}, we introduce the quantity
\[
	C(r):=v_0^2(r)+\frac{2\l m_0(r)}{(n-2)r^{n-2}}
\]
for $n \ge 3$.
It can be said that this describes the balance between the initial velocity and the
strength of the force governed by the Poisson equation.
This clarifies the description of the conditions for global existence.
The large time behavior of $X$ is also distinguished by $C$ (Remark \ref{rmk:orderX}).
For the proof of Theorem \ref{thm:EP}, we use two propositions
(Propositions \ref{prop:attractive} and
\ref{prop:3Dl-}, below).
We first prove Proposition \ref{prop:attractive} and then prove the theorem.

\subsection{The necessary and sufficient condition for the attractive case}
We first consider the case $\l<0$.
We use the function space $D^k$ defined in \eqref{def:Dk}.
The following result is announced but not proven in \cite{MaRIMS}.
\begin{proposition}[Critical thresholds for $\l<0$ case]\label{prop:attractive}
Suppose $\l<0$, $n\ge1$, $\rho_0\in D^k$, and $v_0\in D^{k+1}$ with $v_0(0)=0$
for an integer $k\ge0$.
\begin{enumerate}
\item If $n =1$ or $2$ then the solution to \eqref{eq:a:rEP} is global if and only if
$\rho_0 (r)= 0$, $v_0(r)\ge0$, and $v_0^\prime(r) \ge 0$ hold
for all $r \ge 0$.
In particular, if $\rho_0 \not\equiv 0$ then the solution breaks down in finite time.
\item If $n \ge3$ then the solution is global if and only if
$v_0(r) \ge 0$, $C(r) \ge 0$, and $C^\prime(r) \ge 0$ hold for all $r \ge0$.
\end{enumerate}
If $\rho_0$ and $v_0$ satisfy the condition for global existence, then
the corresponding solution of \eqref{eq:a:rEP} satisfies
\begin{align*}
	\rho &{}\in C^2([0,\I),D^m) \cap C^\I((0,\I),D^m), \\
	v &{}\in C^1([0,\I),D^{m+1}) \cap C^\I((0,\I),D^{m+1}).
\end{align*}
The solution is unique in $C^2([0,\I),D^0) \times C^2([0,\I),D^1)$
and also solves \eqref{eq:a:EP} in the distribution sense.
\end{proposition}
\begin{proof}
Let us recall some facts from \cite{MaRIMS}.
Under the assumptions of Proposition \ref{prop:attractive},
we deduce from Proposition 2.4 in \cite{MaRIMS} that
\eqref{eq:a:rEP} has a unique solution
\begin{align*}
	\rho &{}\in C^2([0,T),D^k) \cap C^\I((0,T),D^k), \\
	v &{}\in C^1([0,T),D^{k+1}) \cap C^\I((0,T),D^{k+1}),
\end{align*}
provided there exists $T$ such that, for $t\in[0,T]$,
$X(t,R)>0$, $R>0$ and $B(t,R)>0$, $R\ge0$.
Moreover, if $B(t_c,R_c)=0$ holds for some $(t_c,R_c)$ then
the solution breaks down at $t=t_c$ (see, Corollary 5.2 in \cite{ELT-IUMJ} or Proposition 2.3 in \cite{MaRIMS}).
Furthermore, if $X(t_0,R_0)=0$ for some $(t_0,R_0)$ then
such $(t_c, R_c)$ exists and $t_c\le t_0$, $R_c\le R_0$ (see, Lemma 2.9 in \cite{MaRIMS}).

{\bf Step 1.}
We begin with the one-dimensional case.
If $\rho_0$ is not identically zero, then we can choose $R_0$ so that $m_0(R_0)>0$.
Twice integration of \eqref{eq:ODE1} yields
$X(t,R_0)=R_0+v_0(R_0)t-(|\l| m_0(R_0)/2)t^2$.
Therefore, we can find $t_0$ such that $X(t_0,R_0)=0$, 
which leads to the finite-time breakdown of the solution.
On the other hand, if $\rho_0\equiv0$ then $X(t,R)=R+v_0(R)t$ and $B(t,R)=1+v_0^\prime(R)t$.
Hence, the solution is global if and only if $v_0(R) \ge0$
and $v_0^\prime(R)\ge0$ holds for all $R\ge0$.

{\bf Step 2.}
We next treat the two-dimensional case.
If there exists $R_0$ such that $m_0(R_0)>0$.
Multiplying \eqref{eq:ODE1} by $X^\prime$, we obtain
\[
	0 \le (X^\prime(t,R_0))^2=v_0(R_0)^2 - 2|\l| m_0(R_0) \log \(\frac{X(t,R_0)}{R_0}\),
\]
which yields an upper bound of $X$:
\[
	X(t,R_0) \le R_0 \exp \(\frac{v_0(R_0)^2}{2|\l| m_0(R_0)}\) =: X_{\mathrm{ub}}.
\]
Plugging this to \eqref{eq:ODE1}, we see that
\[
	X^{\prime\prime}(t,R_0) \le -\frac{|\l| m_0(R_0)}{ X_{\mathrm{ub}}} <0.
\]
Therefore, there exists $t_0$ such that $X(t_0,R_0)=0$.
In the case where $\rho\equiv0$, by the same argument as in the one-dimensional case,
we see that the solution is global if and only if 
$v_0(R)\ge0$ and $v_0^\prime(R)\ge0$ hold for all $R\ge0$.

{\bf Step 3.}
Let us proceed to $n\ge 3$ case.
For simplicity, we use
\[
	A(r):=\frac{2|\l| m_0(r)}{n-2}\ge 0.
\]
Notice that $C$ is written as
$C(r)=v_0(r)^2 -{A(r)}/{r^{n-2}}$.
We first note that $v_0\ge0$ is necessary for global existence.
Indeed, if $v_0(R_0)<0$ for some $R_0>0$,
then $X^{\prime\prime}(t,R_0)\le0$ follows from \eqref{eq:ODE1} and so
$X^\prime(t,R) \le X^\prime(0,R)=v_0(R)<0$.
Hence, there exists $t_0$ such that $X(t_0,R_0)=0$.
We next show that $C \ge 0$ is also necessary for global existence.
Assume that there exists $R_0$ such that $C(R_0)<0$.
In this case, $A(R_0)>0$ by definition of $C$.
Then, multiplying \eqref{eq:ODE1} by $X^\prime$, we obtain
\[
	0\le (X^\prime(t,R_0))^2 = C(R_0) + \frac{A(R_0)}{X(t,R_0)^{n-2}}.
\]
This yields an upper bound of $X$:
\[
	X(t,R_0) \le \(\frac{|C(R_0)|}{A(R_0)}\)^{\frac1{n-2}}.
\]
Then, the same argument as in the two-dimensional case shows 
the existence of $t_0$ such that $X(t_0,R_0)=0$.
Therefore, $C \ge 0$ is necessary for global existence.

In the following, we suppose $v_0\ge0$ and $C\ge0$ are satisfied.
Under this restriction, let us show that the solution is global if and only if
$C^\prime(R)\ge0$ holds for all $R\ge0$.
What to show is that
\begin{equation}\label{eq:PFBattractive}
	C^\prime(R) \ge 0 \Longleftrightarrow B(t,R)>0, \quad \forall t \ge0.
\end{equation}
We first consider the case $v_0(R)>0$.
Then, $C(R)>0$ or $A(R)>0$ hold.
Moreover, $X(t,R)\to \I$ as $t\to\I$ since $X^{\prime\prime}(t,R) \ge0$
and so $X^\prime(t,R)\ge X^\prime(0,R)=v_0(R)>0$.
In this case, by multiplication of \eqref{eq:ODE1} with $X^\prime$,
\[
	X^\prime(t,R) = \sqrt{C(R) + \frac{A(R)}{X(t,R)^{n-2}}}>0,
\]
and so
\[
	\int_R^{X(t,R)} \frac{dy}{\sqrt{C(R) + A(R)y^{-(n-2)}}}=t.
\]
Differentiate with respect to $R$ to obtain
\[
	\frac{B(t,R)}{X^\prime(t,R)}
	-\frac{1}{v_0(R)}
	-\frac12 \int_R^{X(t,R)} \frac{C^\prime(R)
	+ A^\prime(R) y^{-(n-2)}}{\(C(R) + A(R)y^{-(n-2)}\)^{3/2}}dy=0.
\]
We put
\[
	\widetilde{B}(t,R):=\frac{B(t,R)}{X^\prime(t,R)}=
	\frac{1}{v_0(R)}
	+\frac12 \int_R^{X(t,R)} \frac{C^\prime(R)
	+ A^\prime(R) y^{-(n-2)}}{\(C(R) + A(R)y^{-(n-2)}\)^{3/2}}dy.
\]
Two quantity $\widetilde{B}$ and $B$ have the same sign.
Notice that 
\[
	A^\prime(R)=\frac{2\l}{n-2}\rho_0(R)R^{n-1}\ge 0
\]
and that the denominator in the last integral is always positive.
Therefore, if $C^\prime (R)\ge 0$ then the above integral is nonnegative,
and so $\widetilde{B}(t,R)$ stays positive for all $t\ge0$.
On the other hand, if $C^\prime(R) < 0$ then the integral in $\widetilde{B}(t,R)$
tends to $-\I$ as $t\to\I$.
This is because, 
choosing $X_0$ so large that $C^\prime(R)+ A^\prime(R)X_0^{-(n-2)}<-|C^\prime(R)|/2$,
we have
\[
	\int_{X_0}^{X(t,R)} \frac{C^\prime(R)
	+ A^\prime(R) y^{-(n-2)}}{\(C(R) + A(R)y^{-(n-2)}\)^{3/2}}dy
	< -\int_{X_0}^{X(t,R)} \frac{|C^\prime(R)|y^{3(n-2)/2}}{2A(R)^{3/2}} dy
\]
if $A(R)>0$ and
\[
	\int_{X_0}^{X(t,R)} \frac{C^\prime(R)
	+ A^\prime(R) y^{-(n-2)}}{\(C(R) + A(R)y^{-(n-2)}\)^{3/2}}dy
	< -\int_{X_0}^{X(t,R)} \frac{|C^\prime(R)|}{2C(R)^{3/2}} dy
\]
if $C(R)>0$.
The right hand sides of both inequalities tend to $-\I$ as $t\to\I$.
Therefore, we can choose $t_c$ such that $B(t_c,R)=0$.

We finally discuss the case where $v_0(R)=0$.
In this case, since $C(R)\ge0$, 
we have $C(R)=0$ and so $A(R)=0$ ($m_0(R)=0$) by the definition of $C$.
It implies that $\rho(r) = 0$ for all $r \le R$ and
so that, for all $r\le R$, $X^\prime(t,r)\equiv 0$ and $X(t,r)\equiv r$.
Hence, by continuity of $B$, one verifies that $B(t,R)=\lim_{r\uparrow R} \partial_R X(t,r) = 1>0$
for all $t\ge0$. Note that $C^\prime(R)=0$ since $\rho(R)=0$.
Thus, \eqref{eq:PFBattractive} is justified.
\end{proof}

\subsection{The necessary condition for the repulsive case}
The following result is Remark 5.4 of \cite{ELT-IUMJ} if we restrict our attention to 
the case where $v_0 \ge 0$, and this is also a part of Theorem 1.7 in \cite{MaRIMS}.
\begin{proposition}[Necessary condition for $\l>0$ case]\label{prop:3Dl-}
Let $\l>0$, $n\ge3$, $\rho_0\in D^0$, and $v_0\in D^1$ with $v_0(0)=0$.
Then, the classical solution of \eqref{eq:a:rEP} is global
only if $C^\prime(R)\ge0$ for all $R\ge0$.
\end{proposition}
\subsection{Proof of Theorem \ref{thm:EP}}

\begin{proof}[Proof of Theorem \ref{thm:EP}]
By Proposition \ref{prop:attractive}, the solution breaks down in finite time
if $n=1$, $2$ and $\l<0$, since $\rho_0\not\equiv0$.
Suppose $n\ge 3$.
By assumptions on the initial data, we have $C(0)=0$ and $C(r)\to 0$ as $r\to\I$.
Now, Propositions \ref{prop:attractive} and \ref{prop:3Dl-} imply that the solution is global only if $C^\prime(R) \ge 0$ for all $R\ge0$, that is,
only if $C \equiv 0$.
In the $\l<0$ case, Proposition \ref{prop:attractive} shows
the solution is global if we take the positive root:
\begin{equation}\label{eq:globalv}
	v_0(R) = \sqrt{\frac{A(R)}{R^{n-2}}} \ge 0.
\end{equation}
If $\l>0$ then $C \equiv 0$ implies $\rho\equiv0$, which is excluded by assumption.

If $n\ge 3$, $\l<0$, and $v_0(R)$ is given as \eqref{eq:globalv} then
$C\equiv0$ and so $X$ satisfies the equation
\[
	X^\prime(t,R)= \sqrt{\frac{A(R)}{X(t,R)^{n-2}}}, \quad X(0,R)=R.
\]
By separation of variables, we obtain
$X(t,R) = R(1 + \frac{nv_0(R)}{2R}t)^{2/n}$.
Then, \eqref{eq:Xtorho} gives the explicit representation of the solution.
\end{proof}
\begin{remark}\label{rmk:orderX}
The value $C(R)$ is useful to
describe the large time behavior of $X(t,R)$ for $n\ge 3$.
In previous results, we have already established the estimate
\[
	C_1 t^{2/n}+o(t^{2/n}) \le X(t,R) \le C_2 t + o(t)
\]
for a constant $C_1$ and $C_2=\sqrt{C(R)}$ in some cases
(see \cite[Remark 5.1]{ELT-IUMJ}).
Notice that the lower bound is $O(t^{2/n})$ and the upper bound is $O(t)$ as $t\to\I$.
This estimate is sharp in such a sense that, as $t\to\I$,
the both cases $X(t,R)=O(t^{2/n})$ and $X(t,R)=O(t)$ can happen.
Now, we summarize as follows:
Let $n\ge 3$, $\l \in \R\setminus\{0\}$, $v_0(R) \in \R$, $m_0(R)> 0$, and
\[
	C(R)=v_0(R)^2 + \frac{2\l m_0(R)}{(n-2)R^{n-2}}.
\]
Let $X(t,R)$ be a solution of \eqref{eq:ODE1}.
\begin{itemize}
\item If $C(R)>0$, then $X(t,R)>0$ for all $t\ge 0$ and $X(t,R)=\sqrt{C(R)}t + o(t)$
as $t\to\I$.
\item If $C(R)=0$, then $\l<0$ and
\[
	X(t,R) = R\( 1+ t\sqrt{\frac{|\l |n^2m_0(R)}{2(n-2)R^n}} \)^{\frac2n}=O(t^{\frac2n})
\]
as $t\to\I$.
\item If $C(R)<0$, then $\l<0$ and there exists $t_c<\I$ such that $X(t_c,R)=0$.
\end{itemize}
\end{remark}

\subsection*{Acknowledgments}
The author expresses his deep gratitude
to Professor Remi Carles for fruitful discussions in Kyoto.
Deep appreciation goes to Professors Yoshio Tsutsumi and Hideo Kubo
for their valuable advice and constant encouragement.
This research is supported by JSPS fellow.

\bibliographystyle{amsplain}
\bibliography{caustic}

\end{document}